\documentclass[final,leqno,letterpaper]{etna}

\setbibdata{1}{xx}{46}{2017} 

\usepackage{graphicx}%
\usepackage{amsmath,amssymb,amsfonts}%
\usepackage{mathrsfs}%
\usepackage{xcolor}%
\usepackage{textcomp}%
\usepackage{manyfoot}%
\usepackage{booktabs}%
\usepackage{algorithm}%
\usepackage{algorithmicx}%
\usepackage{algpseudocode}%
\usepackage{listings}%

\let\Re\relax
\DeclareMathOperator{\Re}{\mathrm{Re}}
\let\Im\relax
\DeclareMathOperator{\Im}{\mathrm{Im}}

\renewcommand{\i}{\mathrm{i}}
\newcommand{\bI}{{\bf I}}
\newcommand{\bM}{{\bf M}}
\newcommand{\bN}{{\bf N}}
\newcommand{\CC}{{\mathbb C}}
\newcommand{\RR}{{\mathbb R}}

\newcommand{\HH}{{\mathbb H}}

\usepackage{lineno}

\hypersetup{%
	pdftitle={A computational method for multiple steady Hele-Shaw bubbles in  planar domains},
	pdfauthor={Mohamed M.S. Nasser, Christopher C. Green, El Mostafa Kalmoun},
	pdfkeywords={free boundary problems, Hele-Shaw flows, generalized Neumann kernel, multiply connected domains}
}

\title{A computational method for multiple steady Hele-Shaw bubbles in  planar domains\thanks{%
		Received... Accepted... Published online on... Recommended by....
}}

\author{Mohamed M.S. Nasser\footnotemark[2]
	\and Christopher C. Green\footnotemark[2]
	\and El Mostafa Kalmoun\footnotemark[3]}

\shorttitle{Hele-Shaw bubbles in planar domains} 
\shortauthor{M.~Nasser, C.~Green \& E. Kalmoun}

\begin{document}
	
	\maketitle
	
	\renewcommand{\thefootnote}{\fnsymbol{footnote}}
	
	\footnotetext[2]{Department of Mathematics, Statistics \& Physics, Wichita State University, Wichita, KS 67260-0033, USA (mms.nasser@wichita.edu, christopher.green@wichita.edu).}
	\footnotetext[3]{School of Science \& Engineering, Al Akhawayn University in Ifrane, PO Box 104, Ifrane, 53000, Morocco (e.kalmoun@aui.ma).}

\begin{abstract}
We present a unified numerical method to determine the shapes of multiple Hele-Shaw bubbles in steady motion, and in the absence of surface tension, in three planar domains: free space, the upper half-plane, and an infinite channel. Our approach is based on solving the free boundary problem for the bubble boundaries using a fast and accurate boundary integral method. The main advantage of our method is that it allows for the treatment of a very high number of bubbles.
The presented method is validated by recovering some existing results for steady bubbles in channels and free space. Several numerical examples are presented, many of which feature configurations of bubbles that have not appeared in the literature before.
\end{abstract}
	
	

\section{Introduction}

Hele-Shaw flows have provided a major source of interest for mathematicians and physicists over the past century since many physical processes involving moving boundaries can be reduced, after various idealizations, to such flow systems. Indeed, the Hele-Shaw system is a  classical paradigm for Laplacian growth processes \cite{gust,gustV}. In the laboratory, the Hele-Shaw cell apparatus consists of sandwiching viscous fluid between two close-to-touching parallel plates. This configuration enables the study of quasi-two-dimensional fluid flows.  
In industry and in nature, Hele-Shaw flows arise in numerous contexts, from dendritic crystal growth and oil reservoir engineering to pattern formation and bacterial growth in biological systems.
Furthermore, flow in domains such as the upper-half plane and an infinite channel is highly relevant in both engineering and natural settings, where the interaction of fluid flow with boundaries plays a critical role.

In the current work, we focus on Hele-Shaw bubbles steadily translating in three different planar geometries: free space, the upper half-plane, and an infinite horizontal channel. We integrate several different methods in the literature by presenting a fast unified computational method for multiple Hele-Shaw bubbles in steady motion in the aforementioned three geometries. We focus primarily on the solution of the free boundary problems for the bubble boundary shapes. 
An important physical assumption that we will make throughout this work is the absence of surface tension on the bubble boundaries. Without surface tension, there is a continuum of possible bubble speeds. This absence contrasts with cases where surface tension is present which leads to a discrete set of bubble speeds and boundary shapes \cite{GLM}. The governing equation for the flow is Laplace's equation, derived from Darcy’s law~\cite{gust,gustV}. By formulating the free boundary problem in terms of complex potentials, we compute the conformal mapping from a multiply connected circular domain onto the region exterior to multiple bubbles in a particular geometry, using a simple Galilean transformation to link two reference frames.

Over the past few decades, the resolution of Hele-Shaw bubble boundary shapes has received considerable attention. Taylor \& Saffman \cite{tay}, Tanveer \cite{tan} and Vasconcelos \cite{vas94,vas01} considered single bubbles and periodic arrays of bubbles in different configurations. Crowdy~\cite{Cro09} generalized these works by deriving an analytical formula in terms of the Schottky-Klein prime function~\cite{Cro20} to calculate multiple bubble boundary shapes in free space.
The prime function method presented in~\cite{Cro09} can be used also to find an explicit formula for bubbles in the upper half-plane by appealing to the up-down symmetry across the real line and placing image bubbles in the lower half-plane.
Of course, this approach requires considering twice as many bubbles and hence more calculations are necessary.  
An example of this approach for a single bubble configuration is given in~\cite{kha}. 
Green \& Vasconcelos~\cite{GV} also used the function theory of the prime function to derive analytical formulae for multiple bubbles in channel domains. 
Similar analytic formulae could of course be derived using the prime function for other geometries. While the prime function provides an elegant mathematical framework, determining the parameters in the associated analytic formulae requires solving nonlinear systems which become increasingly challenging as the connectivity of the domain increases. 
For the numerical computation of the prime function, the reader is referred to~\cite{CKGN}. Recently, Vasconcelos~\cite{vas15} used secondary prime functions to derive steady solutions describing multiple Hele-Shaw bubbles and fingers in a channel geometry.

In this paper, we employ a well-established boundary integral equation method to calculate pairs of complex potentials describing the bubble flow regimes which affords us accuracy, speed, and the ability to deal with multiply connected domains with a large number of bubbles.
Our proposed method can be used for the three principal planar multiply connected domains we have chosen to consider with only minor modifications needing to be made as we switch between these domains.
For the free space and the infinite channel cases, in addition to the recovery of some of the figures presented in \cite{Cro09,GV,vas15}, we will also consider configurations with a large number of bubbles. 

This paper has the following layout. In section 2, we formulate the free boundary problems governing the shapes of the bubble boundaries using appropriate physical laws. Section 3 reviews the numerical conformal mapping method that we will employ. We then consider several examples of bubble configurations in section 4, each having a different domain connectivity. To illustrate the effectiveness of our method, we adapt our numerical procedure to recover several existing results for Hele-Shaw bubbles in free space~\cite{Cro09} and in the channel~\cite{GV}. Concluding remarks are made in section~5.   
	
\section{Problem formulations} \label{sec:prob}

In this section, we review the formulation of the free boundary problem describing the Hele-Shaw bubble flow in three geometric configurations: free space with $m+1$ bubbles; the upper half-plane with $m$ bubbles; and an infinite channel with $m$ bubbles,
where $m\ge1$. These free boundary problems have been formulated in~\cite{Cro09,kha} for the free space case and in~\cite{Cro09b,GV,tan,tay,vas01,vas15} for the infinite channel case.

Consider viscous fluid flow in one of the three geometric configurations mentioned above and suppose that there are $m$ Hele-Shaw bubbles in steady motion in this flow, where the shapes of these bubble boundaries are unknown \textit{a priori} and must be determined as part of the solution to a governing free boundary problem. Under certain physical assumptions about the bubbles and the flow field (discussed below), our ultimate goal is to solve for the free boundary shapes of the bubbles. We achieve this by formulating the governing free boundary problem in terms of a single complex variable and incorporate the appropriate physics by constructing a pair of complex potential functions in two frames of reference. These complex potentials, when linked together, then provide a conformal mapping formula which can be used to reveal the shapes of the bubble boundaries. 

The motion of the viscous fluid in our Hele-Shaw domain is governed by Darcy's law
\begin{equation}
	{\bf v}=\nabla\phi,
\end{equation}
which asserts that the flow velocity potential $\phi$ is directly proportional to the viscous fluid pressure $p$~\cite{gustV}.
Assuming fluid incompressibility, $\nabla\cdot{\bf v}=0$, it follows at once that $\nabla^2\phi=0$ and hence $\phi$ is a harmonic function. It is therefore natural to formulate the free boundary problem for the bubbles boundaries to be solved in a complex $z$-plane.

We introduce the flow complex potential, $w(z)$, defined by
\[
w(z)=\phi(x,y)+\i\psi(x,y),
\]
where $\psi$ is the flow streamfunction, in a fixed frame of reference, and $z=x+\i y$. If we suppose that the viscous fluid is translating horizontally from left to right with unit speed $V=1$ in the far field, it follows that
\begin{equation}\label{eq:w-inf}
	w(z)\approx  z \quad{\rm as}\quad z\to \infty.
\end{equation}
The complex potential $w(z)$ is analytic everywhere in the viscous fluid domain $D_z$. 

Next, let $\tau(z)$ represent the complex potential in a frame of reference co-traveling with the bubbles. We assume the bubbles are traveling parallel to the $x$-axis at  a constant speed~$U>1$. In this moving frame, $\tau(z)$ is related to $w(z)$ via~\cite[p.~393]{Cro20} (see also~\cite{Cro09,GV}) 
\begin{equation}\label{eq:TW-1}
	\tau(z)=w(z)-U z.
\end{equation}
This Galilean transformation shifts the reference frame to move with the bubbles at constant speed $U$. Consequently, assuming the far field flow has a normalized unit speed, the behavior of $\tau(z)$ at infinity is given by
\[
\tau(z)\approx (1-U)z \quad{\rm as}\quad z\to \infty.
\]

\subsection{The complex potentials}

\subsubsection{Free space}

\begin{figure}[h!]
	\centering
	{\includegraphics[width=0.32\textwidth]{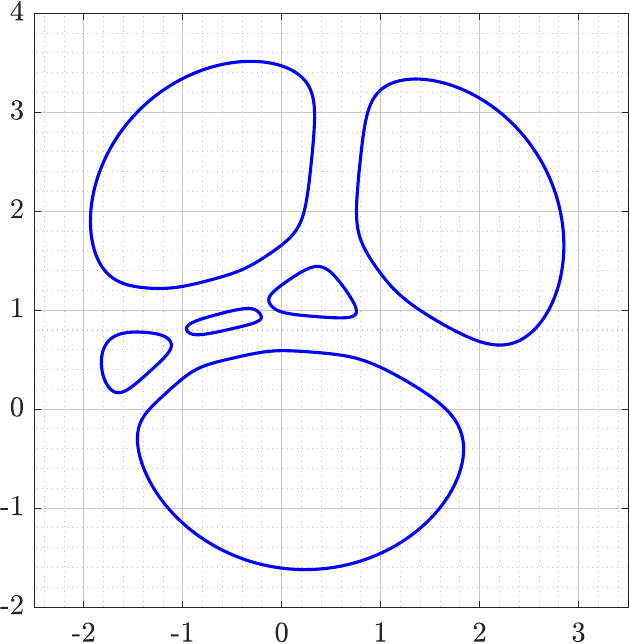}
		\hfill
		\includegraphics[width=0.32\textwidth]{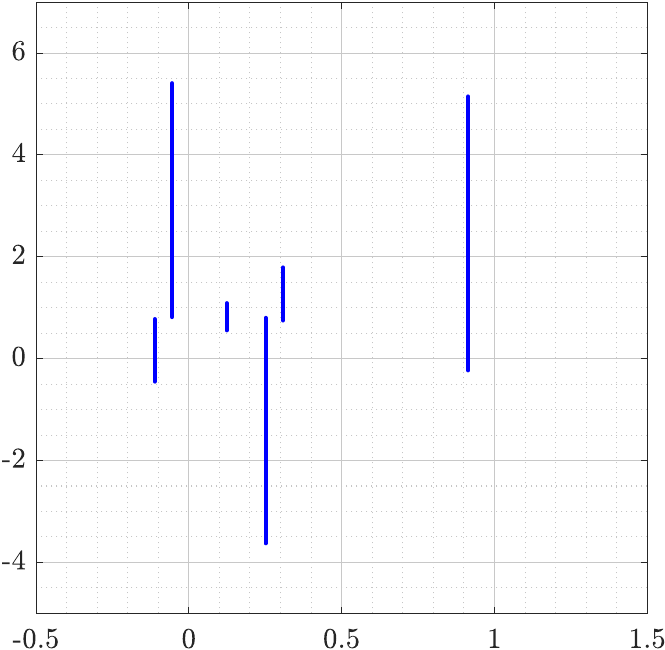}
		\hfill
		\includegraphics[width=0.32\textwidth]{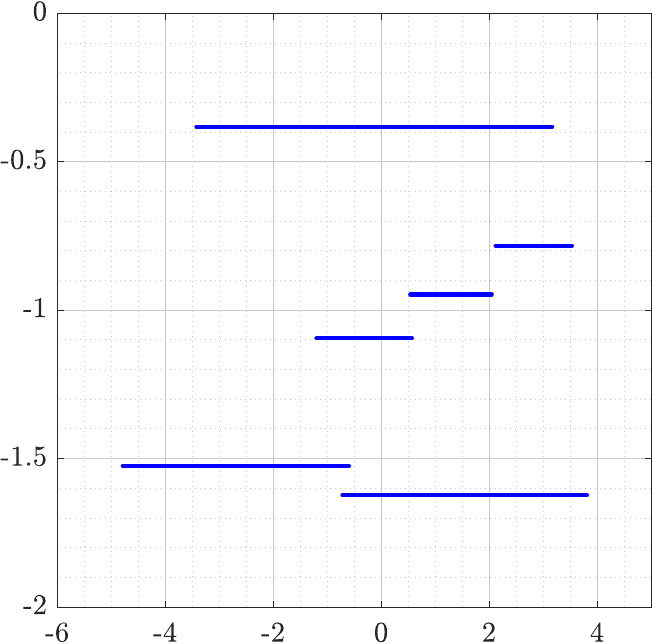}
	}
	\caption{On the left, a schematic of a Hele-Shaw flow in free space containing an assembly of $m+1=6$ bubbles steadily translating with constant speed $U$ parallel to the $x$-axis. The far field fluid speed is assumed to be $V=1$, and the shapes of the bubble boundaries are to be determined. On the middle and on the right, schematic showing the domains in the $w$-plane and in the $\tau$-plane representing the complex potentials in the fixed and moving frames, respectively.}\label{fig:spacefigs}
\end{figure}

The fluid region $D_z$ is the domain exterior to $m+1$ bubbles in the extended complex plane $\CC\cup\{\infty\}$ (see Figure~\ref{fig:spacefigs} (left) for $m=5$). We shall label the boundary of the $j$th bubble by $L_j$, $j=0,1,\ldots,m$.

The complex potential $w(z)$ must satisfy the following boundary conditions~\cite[p.~392]{Cro20}:
\begin{equation}\label{eq:re-cd2-fs}
	\Re[w(z)]={\rm constant} \quad{\rm for}\quad z\in L_j, \quad j=0,1,\ldots,m.
\end{equation}
That is, $\Re[w(z)]$ is constant on each bubble boundary $L_j$ and these constants may be different for each $j$, $j=0,1,\ldots,m$. (This remark applies to all of the boundary conditions in~\eqref{eq:im-cd2-fs}, \eqref{eq:re-cd2}, \eqref{eq:im-cd2}, \eqref{eq:re-cd2-ch}, and~\eqref{eq:im-cd2-ch}, below.)
These conditions~\eqref{eq:re-cd2-fs} follow from the fact that the pressure $p$ must be constant on each bubble boundary.
From~\eqref{eq:re-cd2-fs}, one sees that the flow domain in the $w$-plane is the free space containing $m+1$ vertical rectilinear slits where each slit corresponds to a bubble in the $z$-plane (Figure~\ref{fig:spacefigs} (middle)).

The function $\tau(z)$ satisfies the following boundary conditions~\cite[p.~393]{Cro20}:
\begin{equation}\label{eq:im-cd2-fs}
	\Im[\tau(z)]={\rm constant} \quad{\rm for}\quad z\in L_j, \quad j=1,\ldots,m.
\end{equation}
These conditions~\eqref{eq:im-cd2-fs} mean that, in the co-traveling frame, the bubble boundaries are necessarily streamlines of the flow. It then follows that the flow domain in the $\tau$-plane is the free space containing $m+1$ horizontal rectilinear slits corresponding to the bubbles (Figure~\ref{fig:spacefigs} (right)).

\subsubsection{Upper half-plane}

\begin{figure}[h!]
\centering
{\includegraphics[width=0.32\textwidth]{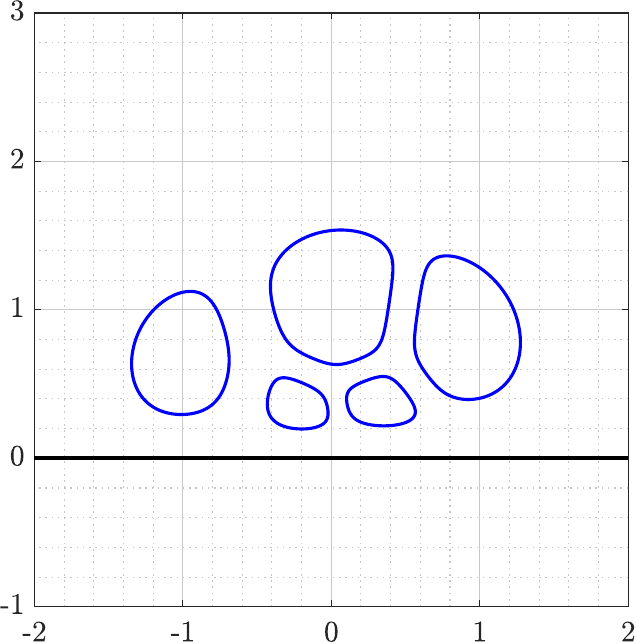}
\hfill
\includegraphics[width=0.32\textwidth]{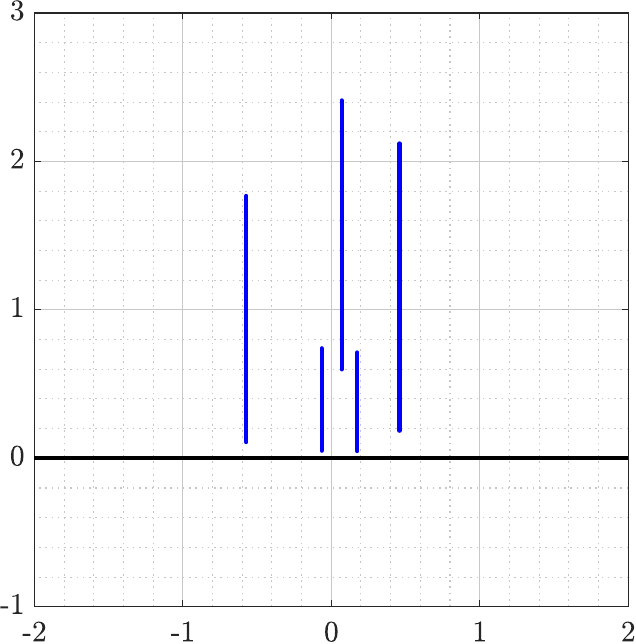}
\hfill
\includegraphics[width=0.32\textwidth]{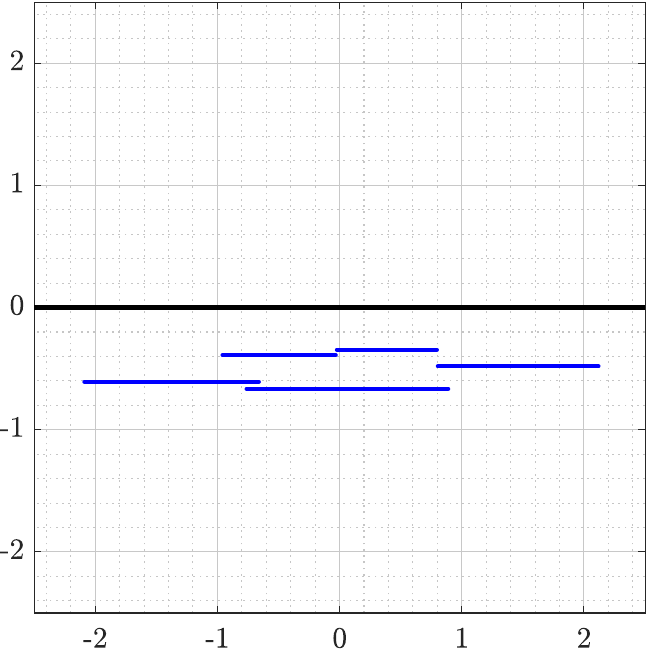}
}
\caption{On the left, a schematic of a Hele-Shaw flow in the upper-half plane containing an assembly of $m=5$ bubbles steadily translating with constant speed $U$ parallel to the $x$-axis. The far field fluid speed is assumed to be $V=1$, and the shapes of the bubble boundaries are to be determined. On the middle and on the right, schematic showing the domains in the $w$-plane and in the $\tau$-plane representing the complex potentials in the fixed and moving frames, respectively.}\label{fig:halffigs}
\end{figure}

The fluid region $D_z$ is the domain exterior to $m$ bubbles in the upper-half plane $\HH=\{z=x+\i y\in\CC\,:\,y>0\}$ (see Figure~\ref{fig:halffigs} (left) when $m=5$). The domain $D_z$ is then bounded by $L=L_0\cup L_1\cup\cdots L_m$ where $L_0=\RR\cup\{\infty\}$.

Similar to the free space case, the complex potential $w(z)$ here must satisfy the following boundary conditions:
\begin{equation}\label{eq:re-cd1}
	\Im[w(z)]=0 \quad{\rm for}\quad z\in \RR,
\end{equation}
and
\begin{equation}\label{eq:re-cd2}
	\Re[w(z)]={\rm constant} \quad{\rm for}\quad z\in L_j, \quad j=1,\ldots,m.
\end{equation}
Condition~\eqref{eq:re-cd1} asserts that the real line $y=0$ is a streamline of the flow.
From conditions~\eqref{eq:re-cd1} and~\eqref{eq:re-cd2}, one sees that the flow domain in the $w$-plane is an upper-half plane containing $m$ vertical rectilinear slits where each slit corresponds to a bubble in the $z$-plane (Figure~\ref{fig:halffigs} (middle)).
Further, the complex potential function $\tau(z)$ in the co-traveling frame with the bubbles  satisfies the following boundary conditions:
\begin{equation}\label{eq:im-cd1}
	\Im[\tau(z)]=0 \quad{\rm for}\quad z\in \RR,
\end{equation}
and
\begin{equation}\label{eq:im-cd2}
	\Im[\tau(z)]={\rm constant} \quad{\rm for}\quad z\in L_j, \quad j=1,\ldots,m.
\end{equation}
Noting that $1-U$ is negative, it then follows from~\eqref{eq:im-cd1} and~\eqref{eq:im-cd2} that the flow domain in the $\tau$-plane is a lower-half plane containing $m$ horizontal rectilinear slits corresponding to the bubbles (Figure~\ref{fig:halffigs} (right)).

\subsubsection{Infinite channel}

\begin{figure}[h!]
\centering
{\includegraphics[width=0.32\textwidth]{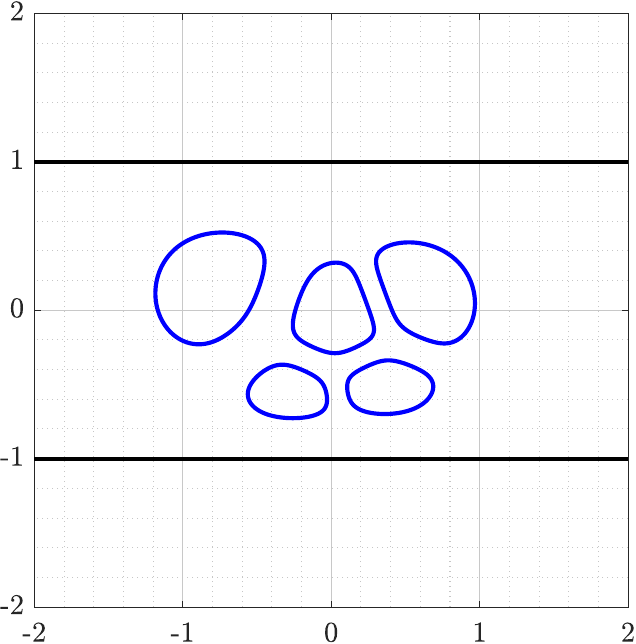}
\hfill
\includegraphics[width=0.32\textwidth]{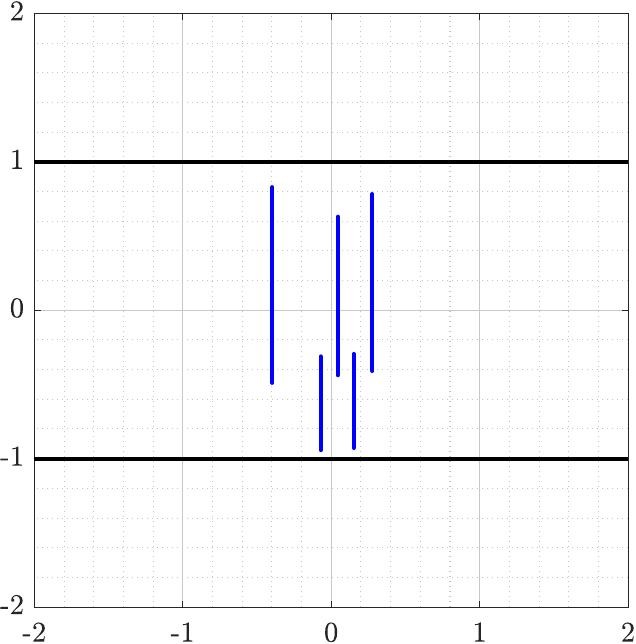}
\hfill
\includegraphics[width=0.32\textwidth]{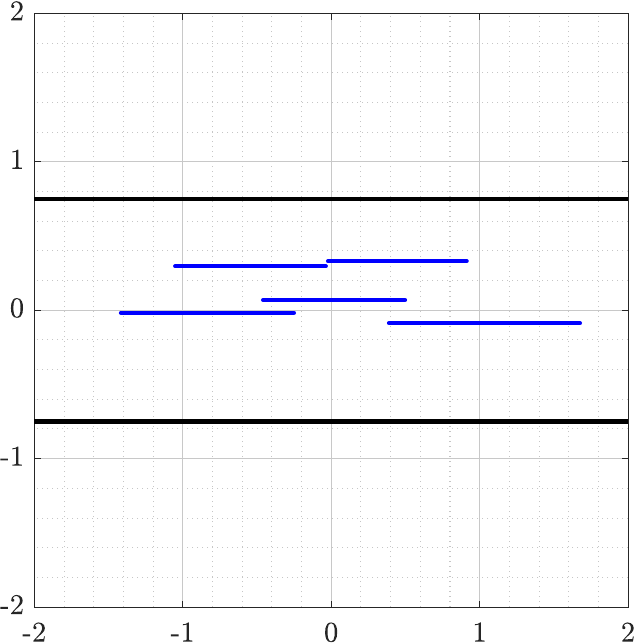}
}
\caption{On the left, a schematic of a Hele-Shaw flow in an infinite channel containing an assembly of $m=5$ bubbles steadily translating with constant speed $U=1.75$ parallel to the $x$-axis. The far field fluid speed is assumed to be $V=1$, and the shapes of the bubble boundaries are to be determined. On the middle and on the right, schematic showing the domains in the $w$-plane and in the $\tau$-plane representing the complex potentials in the fixed and moving frames, respectively.}\label{fig:chanfigs}
\end{figure}

The fluid region $D_z$ is the domain exterior to $m$ bubbles in the infinite channel $S=\{z=x+\i y\in\CC\,:\,-1<y<1\}$ (see Figure~\ref{fig:chanfigs} (left) when $m=5$). The domain $D_z$ is then bounded by $L=L_0\cup L_1\cup\cdots L_m$ where $L_0=\partial S$.

The complex potential $w(z)$ must satisfy the following boundary conditions~\cite{GV}:
\begin{equation}\label{eq:re-cd1-ch}
	\Im[w(z)]=\pm1 \quad{\rm for}\quad \Im z=\pm1,
\end{equation}
and
\begin{equation}\label{eq:re-cd2-ch}
	\Re[w(z)]={\rm constant} \quad{\rm for}\quad z\in L_j, \quad j=1,\ldots,m.
\end{equation}
Conditions~\eqref{eq:re-cd1-ch} assert that the channel walls $\Im z=\pm1$ are streamlines of the flow.
From conditions~\eqref{eq:re-cd1-ch} and~\eqref{eq:re-cd2-ch}, one sees that the flow domain in the $w$-plane is a channel of width $2$ containing $m$ vertical rectilinear slits where each slit corresponds to a bubble in the $z$-plane (Figure~\ref{fig:chanfigs} (middle)).

The function $\tau(z)$ satisfies the following boundary conditions~\cite{GV}:
\begin{equation}\label{eq:im-cd1-ch}
	\Im[\tau(z)]=\pm(1-U) \quad{\rm for}\quad \Im z=\pm1,
\end{equation}
and
\begin{equation}\label{eq:im-cd2-ch}
	\Im[\tau(z)]={\rm constant} \quad{\rm for}\quad z\in L_j, \quad j=1,\ldots,m.
\end{equation}
Since $U>1$, then $1-U<0$ and it follows from~\eqref{eq:im-cd1} and~\eqref{eq:im-cd2} that the flow domain in the $\tau$-plane is a channel of width $2(U-1)$ containing $m$ horizontal rectilinear slits corresponding to the bubbles (see Figure~\ref{fig:chanfigs} (right) for $U=1.75$).

\subsection{The conformal mappings}

\begin{figure}[h!]
	\centering
	{\includegraphics[width=0.32\textwidth]{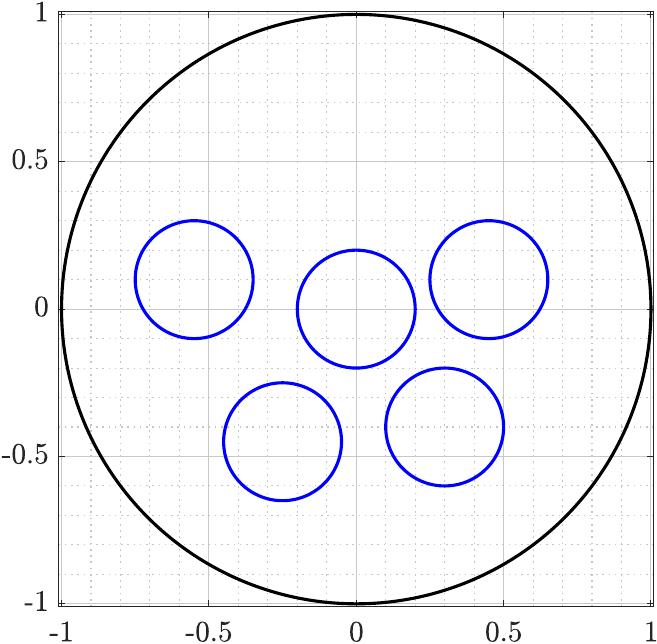}
	}
	\caption{The circular domain $D_\zeta$ in the auxiliary complex $\zeta$-plane.}\label{fig:cirdom}
\end{figure}

Let us now introduce $z(\zeta)$ to be the conformal mapping from a bounded $(m+1)$-connected circular domain $D_\zeta$ in a complex $\zeta$-plane (Figure~\ref{fig:cirdom}) onto the $(m+1)$-connected fluid region $D_z$. 
For the domain $D_\zeta$, we label the unit circle by $C_0$ and the $m$ inner circular boundaries by $C_1,\ldots,C_m$, where the bubble boundary $L_j$ is the image of $C_j$ under the conformal mapping $z(\zeta)$. The center and radius of the circle $C_j$ will be denoted by $z_j$ and $r_j$, respectively.

We define the functions $W(\zeta)$ and $T(\zeta)$ through the following compositions:
\[
W(\zeta)=w(z(\zeta)) \quad{\rm and}\quad T(\zeta)=\tau(z(\zeta)).
\]
It then follows from~\eqref{eq:TW-1} that the conformal mapping $z(\zeta)$ is given by
\begin{equation}\label{eq:Sp-2}
	z(\zeta)=\frac{1}{U}\left(W(\zeta)-T(\zeta)\right).
\end{equation}
The functions $W(\zeta)$ and $T(\zeta)$ are analytic in the circular domain $D_\zeta$ except at a given point $\zeta=\alpha\in D_\zeta$ in the free space case (see~\cite{Cro09,Cro20,GV} for more details).
We have therefore reduced the original free boundary problem in the fluid region $D_z$ to the more manageable task of computing two analytic functions $W(\zeta)$ and $T(\zeta)$ in the circular domain $D_\zeta$ which maps $D_\zeta$ onto $D_z$. 
These functions can be computed using a boundary integral equation method, as we will discuss in the following section.

It is important to set an appropriate length scale of the problem.
Using Green's Theorem, the area $A_j$ of the bubble with boundary $L_j$ is given by
\[
A_j = \Re\left[\frac{1}{2}\oint_{C_j}\i\,\overline{z(\zeta)}z'(\zeta)d\zeta
\right],
\]
where $z(\zeta)$ is the conformal mapping in~\eqref{eq:Sp-2}.
For bubbles in the free space and the upper half-plane domains, note that multiplying the conformal mapping from the circular domain onto the flow domain by a positive real constant will not change the shape of the bubbles, but will affect the scaling and hence the bubble areas. This approach can be used for these two domains to scale the bubbles so that one of these bubbles has a fixed area. For the channel with width $2$, we can not use this approach since doing so would affect the width of the channel.   

\section{Computing the conformal mappings}\label{sec:ie}

Let $D_\zeta$ be the above mentioned bounded multiply connected circular domain of connectivity $m+1$ (see Figure~\ref{fig:cirdom}). In this section, we review a boundary integral equation method from~\cite{Nas-JMAA11,Nas-JMAA13,NF-Siam13} for computing conformal mappings from the circular domain $D_\zeta$ onto each of the following three canonical slit domains: free space with $m+1$ rectilinear slits; the upper half-plane with $m$ rectilinear slits; and an infinite channel with $m$ rectilinear slits.

\subsection{The integral equation}

Let $C_0$ be parametrized by $\zeta_0(t)=e^{\i t}$ for $t\in J_0=[0,2\pi]$, and let $C_j$ be parametrized by $\zeta_j(t)=z_j+r_je^{-\i t}$, $t\in J_j=[0,2\pi]$, for $j=1,2,\ldots,m$.
Let $J$ be the disjoint union of the $m+1$ intervals $J_j=[0,2\pi]$, $j=0,1,\ldots,m$. We define a parametrization of the whole boundary $C=\cup_{j=0}^m C_j$ on $J$ by 
\[
\zeta(t)=\left\{
\begin{array}{cc} 
	\zeta_0(t), & t\in J_0, \\ 
	\zeta_1(t), & t\in J_1, \\
	\vdots \\
	\zeta_m(t), & t\in J_m. \\ 
\end{array}
\right.
\]
Let  $\theta(t)$ also be the piecewise constant real-valued function defined on $J$ by
\begin{equation}\label{eq:thet-1}
	\theta(t)=\theta_j,\quad t\in J_j, \quad j=0,1,\ldots,m,
\end{equation}
where $\theta_j$ is the angle between the slit $S_j$ and the positive real axis, i.e., $\theta_j=0$ for a horizontal slit and $\theta_j=\pi/2$ for a vertical slit. For the upper half-plane and an infinite channel cases, we choose $\theta_0=0$.

With the parametrization $\zeta(t)$ of the whole boundary $C$ and with the above piecewise constant function $\theta(t)$, we define a complex function $A$ by
\begin{equation}\label{eq:A}
	A(t) = 	e^{\i\left(\frac{\pi}{2}-\theta(t)\right)}\,(\zeta(t)-\alpha), \quad t\in J,
\end{equation}
where $\alpha$ is a given point in the domain $D_\zeta$. The generalized Neumann kernel $N(s,t)$ is 
defined for $(s,t)\in J\times J$ by
\begin{equation}\label{eq:N}
	N(s,t) =
	\frac{1}{\pi}\Im\left(\frac{A(s)}{A(t)}\frac{\zeta'(t)}{\zeta(t)-\zeta(s)}\right).
\end{equation}
We also define the following kernel 
\begin{equation}\label{eq:M}
	M(s,t) =
	\frac{1}{\pi}\Re\left(\frac{A(s)}{A(t)}\frac{\zeta'(t)}{\zeta(t)-\zeta(s)}\right),
	\quad (s,t)\in J\times J.
\end{equation}

Let $\nu$ be a fixed number with $0<\nu<1$. In what follows Hölder continuity always means H\"older
continuity with this fixed H\"older exponent $\nu$. 	
Let $H$ be the space of all real functions $\gamma(t)$, $t\in J$, of the form
\[
\gamma(t)=\left\{
\begin{array}{cc} 
	\gamma_0(t), & t\in J_0, \\ 
	\gamma_1(t), & t\in J_1, \\
	\vdots \\
	\gamma_m(t), & t\in J_m, \\ 
\end{array}
\right.
\] 
with real H\"older continuous $2\pi$-periodic functions $\gamma_0,\gamma_1,\ldots,\gamma_m$.
In view of the smoothness of the parametrization $\zeta(t)$ of the boundary $C$, any H\"older continuous real-valued function $\hat\gamma$ on the boundary $C$ can be interpreted via $\gamma(t)=\hat\gamma(\zeta(t))$, $t\in J$, as a function $\gamma\in H$; and vice versa.

The integral operators defined on $H$ with the kernels $N(s,t)$ and $M(s,t)$ are denoted by $\bN$ and $\bM$, respectively. 
The kernel $N(s,t)$ is continuous while the kernel $M(s,t)$ is singular and its singular part is (up to the sign) the Hilbert cotangent kernel~\cite[Lemma~1]{Weg-Nas}. Thus the operator $\bN$ is compact and the operator $\bM$ is singular. Moreover, it follows from~\cite[Lemma~2]{Weg-Nas} that both operators $\bN$ and $\bM$ are bounded in $H$ and both operators map $H$ into $H$ (see the paragraph below the proof of Lemma~2 in~\cite{Weg-Mur-Nas}).
Further details can be found in~\cite{Weg-Mur-Nas,Weg-Nas}.

For a given $\gamma\in H$, there exists a unique solution $\mu\in H$ of the boundary integral equation~\cite[Theorem~2]{Nas-CMFT09}
\begin{equation}\label{eq:ie}
	(\bI-\bN)\mu=-\bM\gamma;
\end{equation}
and a unique function
\begin{equation}\label{eq:h}
	h=[\bM\mu-(\bI-\bN)\gamma]/2
\end{equation}
such that 
\begin{equation}\label{eq:f-half}
	f(\zeta(t))=(\gamma(t)+h(t)+\i\mu(t))/A(t)
\end{equation}
are the boundary values of an analytic function $f(\zeta)$ in the bounded domain $D_\zeta$. The real-valued function $h(t)$ is a piecewise constant function, i.e.,
\[
h(t)=h_j,\quad t\in J_j, \quad j=0,1,\ldots,m,
\]
with real constants $h_0,h_1,\ldots,h_m$. By computing the real functions $\mu$ and $h$, we obtain the boundary values of the analytic function $f(\zeta)$ and then the values of the function $f(\zeta)$ for $\zeta\in D_\zeta$ can be computed by the Cauchy integral formula.

Note that the Fredholm operator $\bN$ and the singular operator $\bM$ are both bounded on $H$ and map $H$ onto itself~\cite{Weg-Nas}. Thus, when $\gamma\in H$, the right-hand side $-\bM\gamma$ of the integral equation~\eqref{eq:ie} is also in $H$. 
Note also that both kernels $M(s,t)$ and $N(s,t)$ depend on the function $A$ which, in turn, depends on the constants $\theta_0,\theta_1,\ldots,\theta_m$. Further, the solvability of the integral equation~\eqref{eq:ie} does not depend on these constants; however, the eigenvalues of the kernel $N(s,t)$ do depend on the values of these constants. See~\cite{Nas-ETNA} for details. 

\subsection{Free space with $m+1$ rectilinear slits}

Let $\Phi(\zeta)$ be the conformal mapping from the circular domain $D_\zeta$ onto the slit domain $\Omega$ consisting of the entire extended complex plane with $m+1$ rectilinear slits such that
\begin{equation}\label{eq:sp-Phi}
\Phi(\alpha)=\infty, \quad \lim_{\zeta\to\alpha}[\Phi(\zeta)-1/(\zeta-\alpha)]=0,
\end{equation}
for a given point $\alpha$ in the circular domain $D_\zeta$.
Each circle $C_j$ is mapped onto a rectilinear slit $S_j$, $j=0,1,\ldots,m$. An example of such a configuration with $m=5$ is shown in Figure~\ref{fig:spacefigs}.

To compute the conformal mapping $\Phi(\zeta)$, we define a function $\gamma(t)$ on $J$ by
\begin{equation}
\gamma(t)=\Im\left[e^{-\i\theta(t)}/(\zeta(t)-\alpha)\right].
\end{equation}
Let $\mu$ be the unique solution of the integral equation~\eqref{eq:ie} and the function $h$ is given by~\eqref{eq:h}, and let $f(\zeta)$ be the analytic function in the domain $D_\zeta$ with the boundary values~\eqref{eq:f-half}. Then the conformal mapping $\Phi(\zeta)$ is given by
\begin{equation}\label{eq:Phi-space}
\Phi(\zeta)=1/(\zeta-\alpha)+(\zeta-\alpha)f(\zeta), \quad \zeta\in D_\zeta\cup C. 
\end{equation}
For more details, we refer the reader to~\cite{Nas-JMAA11}.

By computing the conformal mapping $\Phi(\zeta)$ in~\eqref{eq:Phi-space} with  $\theta_0=\theta_1=\cdots=\theta_m=\pi/2$, the function $W(\zeta)$ in~\eqref{eq:Sp-2} is then given by $W(\zeta)=\Phi(\zeta)$. 
To determine the function $T(\zeta)$ in~\eqref{eq:Sp-2}, we compute the conformal mapping $\Phi(\zeta)$ with $\theta_0=\theta_1=\cdots=\theta_m=0$ and then $T(\zeta)=\Phi(\zeta)$. 
By computing both functions $W(\zeta)$ and $T(\zeta)$, we obtain the conformal mapping $z(\zeta)$ from the circular domain $D_\zeta$ onto the fluid domain $D_z$.

\subsection{Upper half-plane with $m$ rectilinear slits}

Let $\Phi(\zeta)$ be the conformal mapping from the circular domain $D_\zeta$ onto the slit domain $\Omega$ consisting of the upper half-plane with $m$ rectilinear slits such that
\[
z(\i)=\infty.
\]
Refer to Figure~\ref{fig:halffigs} (middle) for an example of such a domain $\Omega$ when $m=5$. 
The external circle $C_0$ will be mapped onto the real line  and each inner circle $C_j$ will be mapped onto a rectilinear slit $S_j$, $j=1,2,\ldots,m$. 

To compute the conformal mapping $\Phi(\zeta)$, we need the following M\"obius transformation
\[
\Psi(\zeta)=\i\frac{\i+\zeta}{\i-\zeta}
\]
which maps the unit circle onto the real line and the interior of the unit circle onto the upper half-plane with $\Psi(\i)=\infty$. We define a function $\gamma(t)$ on $J$ by
\begin{equation}
\gamma(t)=
	\left\{ \begin{array}{l@{\hspace{0.5cm}}l}
		0,&t\in J_{0},\\
		\Im\left[e^{-\i\theta_j}\Psi(\zeta_j(t))\right],&t\in J_j, \quad j=1,2, \ldots,m.
	\end{array}\right.
\end{equation}

Let $\mu$ be the unique solution of the integral equation~\eqref{eq:ie} and the function $h$ is given by~\eqref{eq:h}, and let $f(\zeta)$ be the analytic function in the domain $D_\zeta$ with the boundary values~\eqref{eq:f-half}. Then the conformal mapping $\Phi(\zeta)$ is given by
\begin{equation}\label{eq:Phi-half}
\Phi(\zeta)=\Psi(\zeta)+(\zeta-\alpha)f(\zeta)+\i h_0, \quad \zeta\in G\cup\Gamma. 
\end{equation}
For more details, we refer the reader to~\cite{Nas-JMAA13}.

By computing the conformal mapping $\Phi(\zeta)$ in~\eqref{eq:Phi-half} with $\theta_0=0$ and $\theta_1=\cdots=\theta_m=\pi/2$, the function $W(\zeta)$ in~\eqref{eq:Sp-2} is then given by $W(\zeta)=\Phi(\zeta)$. 
The function $T(\zeta)$ in~\eqref{eq:Sp-2} is determined by computing the conformal mapping $\Phi(\zeta)$ with $\theta_0=\theta_1=\cdots=\theta_m=0$ and then $T(\zeta)=\Phi(\zeta)$. 
Once we obtain both functions $W(\zeta)$ and $T(\zeta)$, the conformal mapping $z(\zeta)$ from the circular domain $D_\zeta$ onto the fluid domain $D_z$ is obtained through~\eqref{eq:Sp-2}.

\subsection{Infinite channel with $m$ rectilinear slits}

Let $\Phi(\zeta)$ be the conformal mapping from the circular domain $D_\zeta$ onto the slit domain $\Omega$ consisting of the infinite channel $|\Im z|<1$ with $m$ rectilinear slits such that
\[
z(\pm1)=\pm\infty+0\i, \quad z(\i)=\i.
\]
Refer to Figure~\ref{fig:chanfigs} for an example of such a domain $\Omega$ when $m=5$. 
The upper half of the external circle $C_0$ will be mapped onto the line $\Im w=1$ and the lower half of the circle $C_0$ will be mapped onto the line $\Im w=-1$. Each inner circle $C_j$ will be mapped onto a rectilinear slit $S_j$, $j=1,2,\ldots,m$. 

To compute the conformal mapping $\Phi(\zeta)$, we need the following transformation
\[
\Psi(\zeta)=\frac{2}{\pi}\log\left(\frac{1+\zeta}{1-\zeta}\right)
\]
which maps the unit disk onto the infinite channel $-1<\Im z<1$ with $\Psi(\pm1)=\pm\infty+0\i$ and $\Psi(\i)=\i$. We define a function $\gamma(t)$ on $J$ by
\begin{equation}
	\gamma(t)=
	\left\{ \begin{array}{l@{\hspace{0.5cm}}l}
		0,&t\in J_{0},\\
		\Im\left[e^{-\i\theta_j}\Psi(\zeta_j(t))\right],&t\in J_j, \quad j=1,2, \ldots,m.
	\end{array}\right.
\end{equation}

Let $\mu$ be the unique solution of the integral equation~\eqref{eq:ie} and the function $h$ is given by\eqref{eq:h}, and let $f(\zeta)$ be the analytic function in the domain $D_\zeta$ with the boundary values~\eqref{eq:f-half}. Then the conformal mapping $\Phi(\zeta)$ is given by
\begin{equation}\label{eq:Phi-chan}
\Phi(\zeta)=\Psi(\zeta)+(\zeta-\alpha)f(\zeta)-(\i-\alpha)f(\i), \quad \zeta\in G\cup\Gamma. 
\end{equation}
For more details, we refer the reader to~\cite{NF-Siam13}.

As in the half-plane case, the function $W(\zeta)$ in~\eqref{eq:Sp-2} can be determined through $W(\zeta)=\Phi(\zeta)$, where $\Phi(\zeta)$ is the conformal mapping in~\eqref{eq:Phi-chan} and is computed with $\theta_0=0$ and $\theta_1=\cdots=\theta_m=\pi/2$.
The function $T(\zeta)$ in~\eqref{eq:Sp-2} is found by computing the conformal mapping $\Phi(\zeta)$ with $\theta_0=\theta_1=\cdots=\theta_m=0$ and then $T(\zeta)=\Phi(\zeta)$.
Once both functions $W(\zeta)$ and $T(\zeta)$ are determined, the conformal mapping $z(\zeta)$ from the circular domain $D_\zeta$ onto the fluid domain $D_z$ is obtained through~\eqref{eq:Sp-2}.

\subsection{The numerical solution of the integral equation}\label{sec:numerical}

The integrands in the integral equation~\eqref{eq:ie} and in the formula~\eqref{eq:h} are $2\pi$-periodic functions.	The smoothness of these integrands depends on the smoothness of the boundary of the
computational domain $D_\zeta$.		
In this paper, the boundary components of the computational domain $D_\zeta$ are circles and the integrands in~\eqref{eq:ie} and~\eqref{eq:h} can be analytically extended to some parallel strip $|\Im t|<\sigma$ in the complex plane. Hence, the trapezoidal rule will then converge exponentially with $\mathcal{O}(e^{-\sigma n})$~\cite{Tre} when it is used to discretize the integrals in~\eqref{eq:ie} and~\eqref{eq:h}.

The function $\bM\gamma$ in the right-hand side of the integral equation~\eqref{eq:ie} is continuous since the function $\gamma$ is H\"older continuous. 
Further, since the integrals in~\eqref{eq:ie} are over $2\pi$-periodic functions, the integral operator $\bN$ can be best discretized by the Nystr\"om method with the trapezoidal rule~\cite{Atk}.
The stability and convergence of the Nystr\"om method is based on the compactness of the operator $\bN$ in the space of continuous functions equipped with the sup-norm, on the convergence of
the trapezoidal rule for all continuous functions, and on the theory of collectively
compact operator sequences (cf.~\cite{Atk}). In this paper, we use the MATLAB function {\tt fbie} from~\cite{Nas-ETNA} to solve numerically the integral equation~\eqref{eq:ie} and to compute the function $h$ in~\eqref{eq:h}.
In {\tt fbie}, the integral equation~\eqref{eq:ie} is discretized by the Nystr\"om method with the trapezoidal rule to obtain an $(m+1)n\times(m+1)n$ linear system where $n$ is the number of the discretization points in each boundary component. 
The numerical solution of the integral equation will then converge with a similar rate of convergence as the trapezoidal rule~\cite[p.~322]{Atk}.	
The linear system obtained by discretizing the integral equation~\eqref{eq:ie} is then solved by the Generalized Minimal Residual (GMRES) method using the MATLAB function {\tt gmres} where the matrix-vector product in {\tt gmres} is computed by the Fast Multipole Method (FMM) via the MATLAB function {\tt zfmm2dpart} from the MATLAB toolbox FMMLIB2D~\cite{Gre-Gim12}. 
The GMRES method is used without restart, with tolerance $10^{-14}$, and with $100$ as the maximal number of allowed iterations. The tolerance of the FMM in {\tt zfmm2dpart} is chosen to be $0.5\times10^{-15}$.
By computing the functions $\mu$ and $h$ numerically, we obtain approximations of the boundary values of the analytic function $f(\zeta)$ through~\eqref{eq:f-half}. The values of $f(\zeta)$ for $\zeta\in D_\zeta$ can be computed by the Cauchy integral formula, and then the values of the conformal mapping $\Phi(\zeta)$ can be computed via~\eqref{eq:Phi-space}, \eqref{eq:Phi-half}, or~\eqref{eq:Phi-chan}.
See~\cite{Nas-ETNA} for details.

\section{Examples}

In this section, we present several examples of Hele-Shaw bubble flow in the three domains of interest, computed using the numerical method described earlier. These examples aim to showcase the versatility and accuracy of our method in handling a wide range of bubble configurations.

\subsection{Free space}

In the first example, we recover the bubbles in an unbounded Hele-Shaw cell found in Crowdy \cite{Cro09}. The parameters of the circular domain $D_\zeta$ used to generate the bubbles in this example are the same as in~\cite[Figures~2--4]{Cro09}.

\begin{example}[Two bubbles]\label{ex:spac2}{\rm
We consider two bubbles which are obtained using the following parameters of the circular domain $D_\zeta$: $r_1=0.4$ and $z_1=0$. When we choose $\alpha=\i\sqrt{r_1}$  in~\eqref{eq:sp-Phi}, the two bubbles will have equal area. The bubbles can then be scaled so that their areas are both $\pi$.		
Figure~\ref{fig:spac2} shows the configurations of the two bubbles for several values of $U$. As the speed $U$ increases, the deformation of the bubbles becomes more pronounced, with the bubbles becoming elongated in the direction of the flow. 
		
In Figure~\ref{fig:spac2-lc}, streamlines are superposed in the frame of reference co-traveling with the bubbles for $U=2$ (left) and $U=4$ (right). 	The pair of stagnation points on each bubble is easily identified by a particular streamline. 			
This example illustrates how the shape of the bubbles and the flow field evolve as the speed $U$ increases. 

\begin{figure}[h!]
	\centering
	{\includegraphics[width=0.4\textwidth]{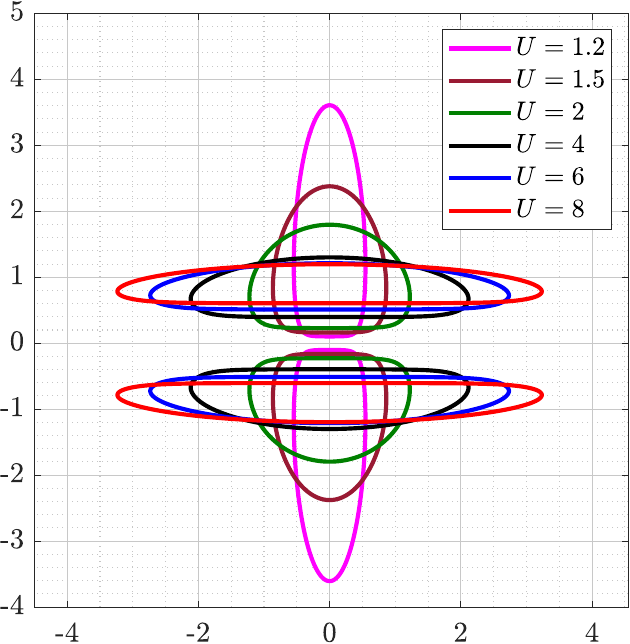}}
	\caption{The configurations of two bubbles in free space in Example~\ref{ex:spac2} corresponding to several values of $U$. The area of each bubble is always $\pi$.}\label{fig:spac2}
\end{figure}

\begin{figure}[h!]
	\centering
	{\hfill\includegraphics[width=0.32\textwidth]{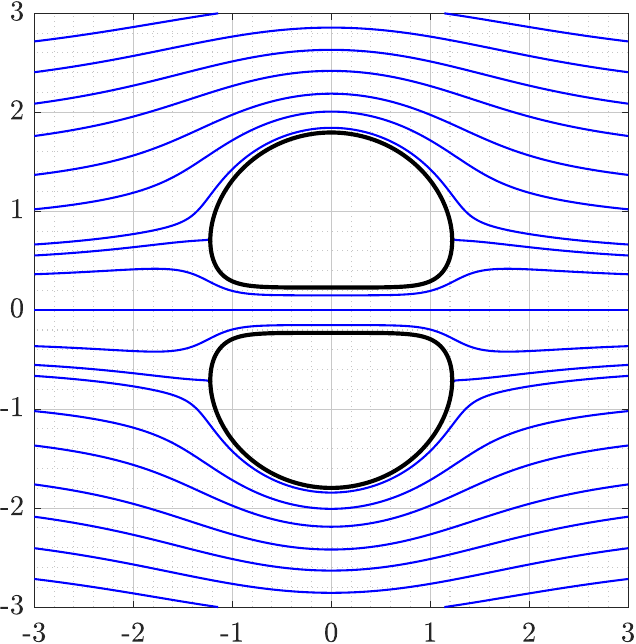}
		\hfill
		\includegraphics[width=0.32\textwidth]{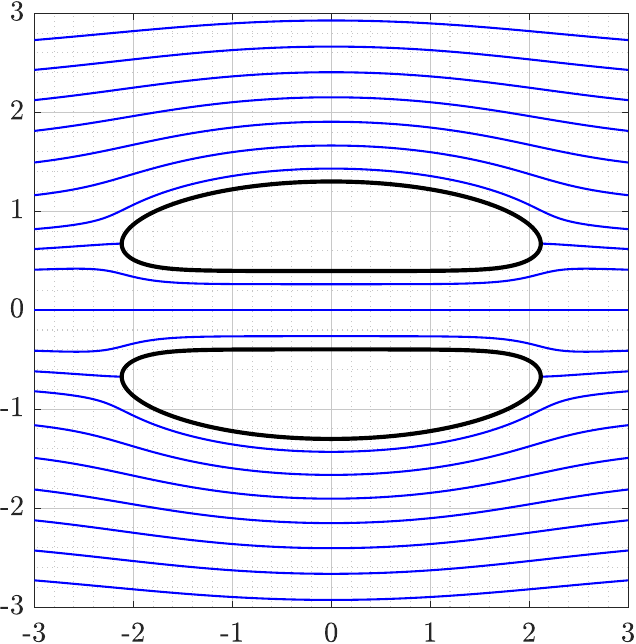}\hfill}
	\caption{The configurations of two bubbles in free space in Example~\ref{ex:spac2} with a set of streamlines superposed in the frame of reference co-traveling with the bubbles when $U=2$ (left) and $U=4$ (right). }\label{fig:spac2-lc}
\end{figure}
}
\end{example}

\begin{example}[Twenty-five bubbles]\label{ex:spac25}{\rm
To illustrate the versatility and effectiveness of our numerical method, we now examine a more complicated setup involving twenty-five bubbles in free space. The parameters of the circular domain $D_\zeta$ used to produce both configurations in Figure~\ref{fig:spacefigLc25} can be found in the MATLAB codes available at \url{https://github.com/mmsnasser/bubbles}. 
		
Figure~\ref{fig:spacefigLc25} presents two configurations of bubbles for speeds $U=2$ (left) and $U=6$ (right), with streamlines superimposed in the frame of reference co-moving with the bubbles. At the higher speed, the bubbles become more horizontally elongated along the flow direction, and the streamline patterns around them simplify. This simplification is likely due to the fluid flow becoming more uniform at higher velocities, with fewer disturbances around the bubbles.
			
These results highlight the ability of the method to handle large assemblies of bubbles. The evolution of bubble deformation and flow patterns with increasing speed is clearly visible. Compared to the previous two-bubble example, the presence of twenty-five asymmetric bubbles introduces significantly more interactions between the bubbles and the surrounding fluid. The varied sizes and positions of the bubbles lead to intricate flow patterns across the bubble array, particularly at lower velocities.

\begin{figure}[h!]
	\centering
	\scalebox{0.43}{\includegraphics[trim=0cm 0.6cm 0cm 1.0cm,clip]{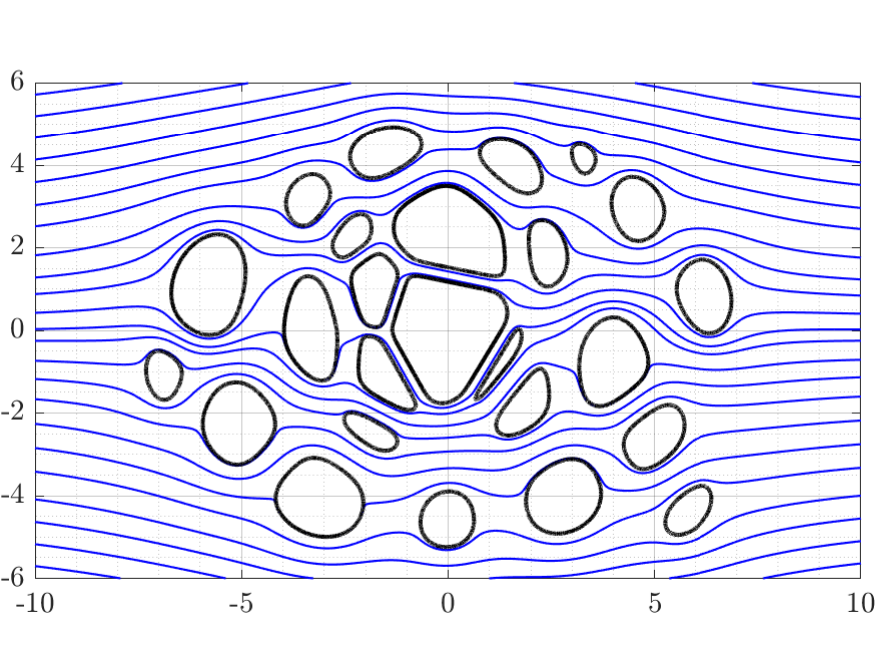}}
	\scalebox{0.43}{\includegraphics[trim=0cm 0.6cm 0cm 1.0cm,clip]{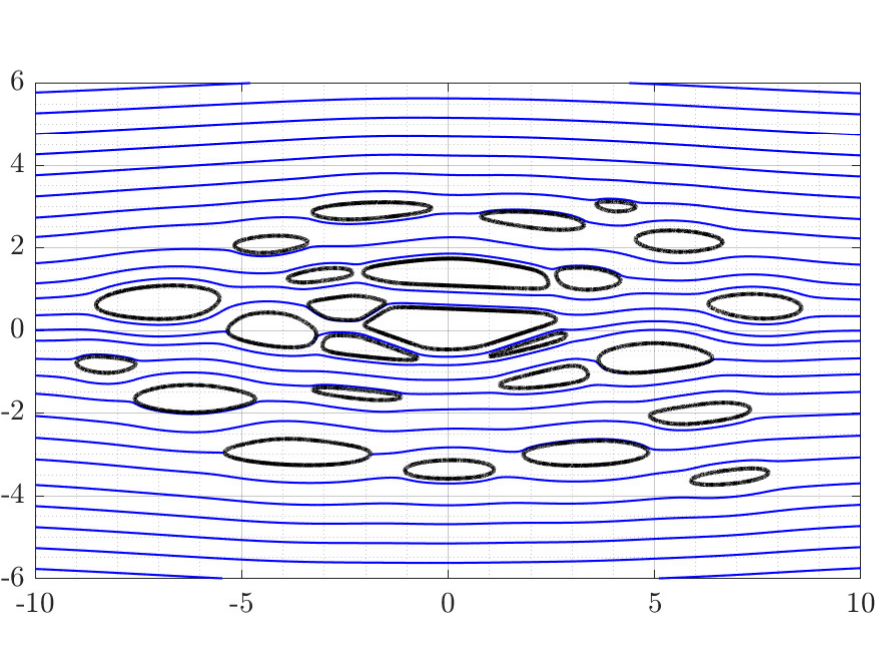}}
	\caption{An example of 25 bubbles in a general asymmetric configuration  in the free space with a set of streamlines superposed in the reference frame co-traveling with the bubbles when $U=2$ (left) and $U=6$ (right).}\label{fig:spacefigLc25}
\end{figure}
}
\end{example}

\subsection{Upper half-plane}

In this subsection, we explore Hele-Shaw bubble flows in the upper half-plane. We consider configurations of bubbles in this region to demonstrate how the geometry of the domain influences the shape and interaction of the bubbles.
	
\begin{figure}[h!]
	\centering
	\scalebox{0.49}{\includegraphics[trim=0cm 0.0cm 0cm 0.0cm,clip]{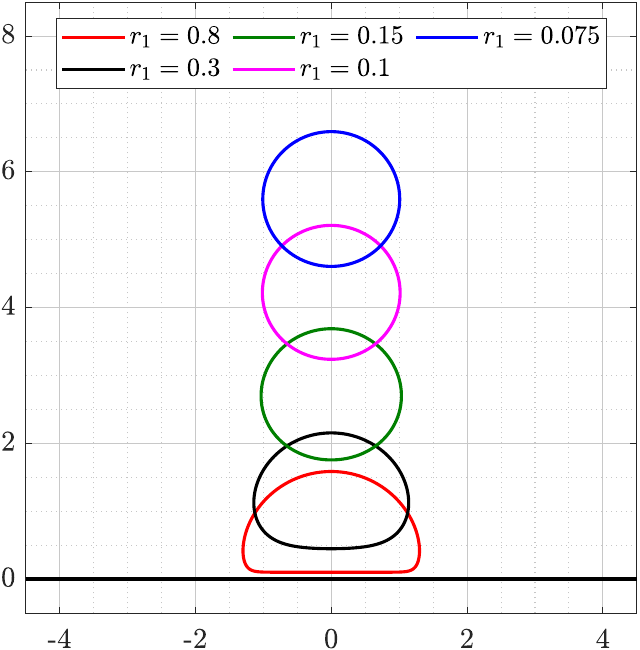}}\hfil
	\scalebox{0.49}{\includegraphics[trim=0cm 0.0cm 0cm 0.0cm,clip]{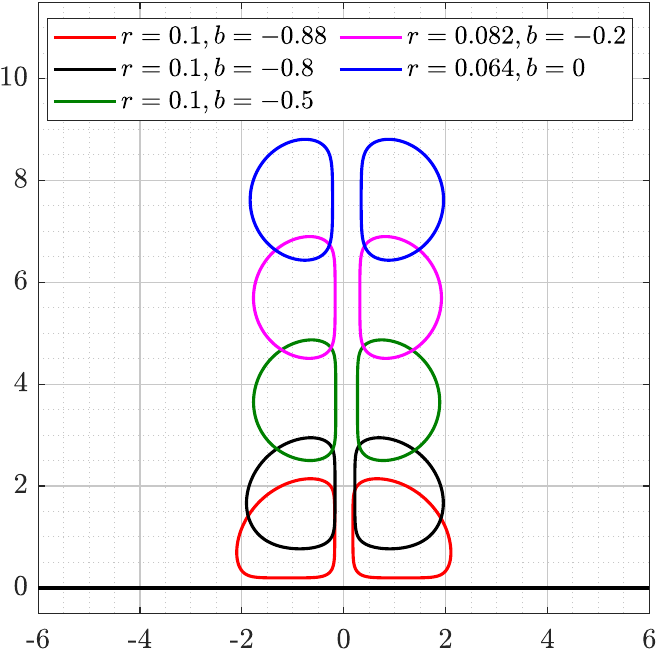}}
	\caption{On the left, the configurations are the single bubbles in the upper half-plane in Example~\ref{ex:1b} corresponding to several values of the parameter $r_1$ when $U=2$. 
		On the right, the configuration are the two bubbles in the upper half-plane in Example~\ref{ex:2bn} corresponding to several values of the parameters $r$ and $b$ when $U=2$. 	
		The area of each of these bubbles is always taken to be $\pi$.}\label{fig:halffig1D}
\end{figure}

\begin{example}[One bubble]\label{ex:1b}{\rm
First we start with an example of a single bubble in the upper half-plane, which constitutes a doubly connected domain. 
In this case, the circular domain $D_\zeta$ is a concentric annulus; $D_\zeta=\{\zeta\,|\,r_1<|\zeta|<1\}$. Here, since we have only one bubble, it can be scaled such that its area is $\pi$.
Figure~\ref{fig:halffig1D} (left) shows the configurations for several values of the parameter $r_1$ when $U=2$, which correspond to different centroid positions above the real line. It is clear that when the bubble is far from the real line, its shape is almost purely circular, gradually becoming semi-circular as it approaches the real line, as expected~\cite{gustV}. 
It is clearly seen that the primary factor influencing the shape of the bubble when it is far from the wall is its interaction with the flow field, rather than boundary effects. In this case, the shape of the bubble tends to approximate that seen in unbounded domains. 
		
With only one bubble in the upper half-plane, we can appeal to the up-down symmetry across the real line when an image bubble is placed in the lower half-plane. As explained in~\cite{kha}, one can use the prime function method presented in~\cite{Cro09} for bubble configurations in free space to derive an explicit formula for this single bubble configuration in the upper half-plane.
}\end{example}

\begin{example}[Two bubbles I]\label{ex:2bn}{\rm
Next, we consider an example of two bubbles in the upper half-plane, forming a triply connected domain. Due to the symmetry, the two bubbles will have the same area and hence both bubbles can be scaled such that their areas are $\pi$.
The parameters of the circular domain $D_\zeta$ are $r_1=r_2=r$ and $z_2=-z_1=(r+0.1)+\i b$. We explore the solutions corresponding to five sets of parameters for $r$ and $b$: $(r,b)=(0.1,-0.88)$, $(r,b)=(0.1,-0.8)$, $(r,b)=(0.1,-0.5)$, $(r,b)=(0.082,-0.2)$, and $(r,b)=(0.064,0)$.
Figure~\ref{fig:halffig1D} (right) shows the configurations for these values of $r$ and $b$ when $U=2$, which correspond to different centroid positions above the real line. It is clear that when the bubbles are far from the real line but close to each other, their shapes are nearly semi-circular. As the bubbles move closer to the real line, their shapes become more quarter-circular.
}\end{example}

\begin{figure}[h!]
	\centering
	{\hfill\includegraphics[width=0.32\textwidth]{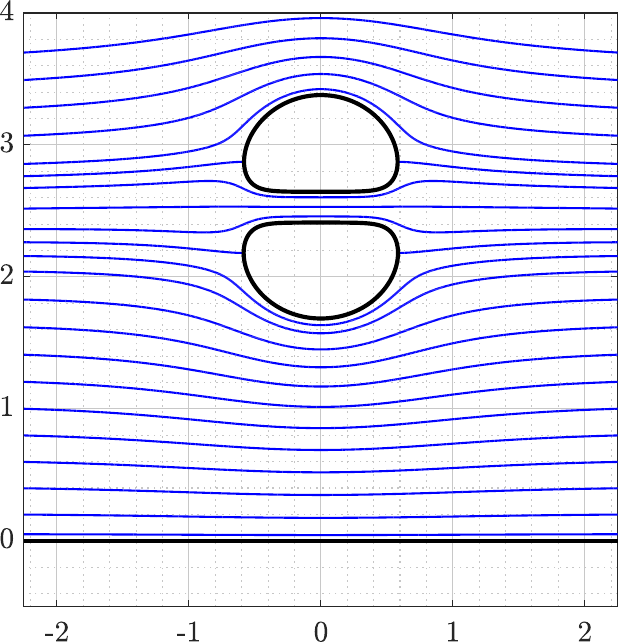}
		\hfill
		\includegraphics[width=0.32\textwidth]{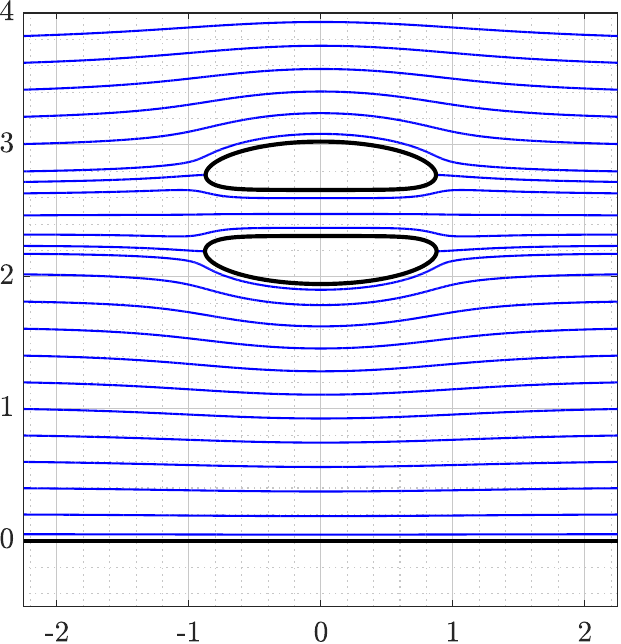}\hfill}
	\caption{The two bubbles in the first case of Example~\ref{ex:2b} with a set of streamlines superposed in the reference frame co-traveling with the bubbles for $U=2$ (left) and $U=4$ (right).}\label{fig:halffig2U}
\end{figure}

\begin{figure}[h!]
	\centering
	{\hfill\includegraphics[width=0.49\textwidth]{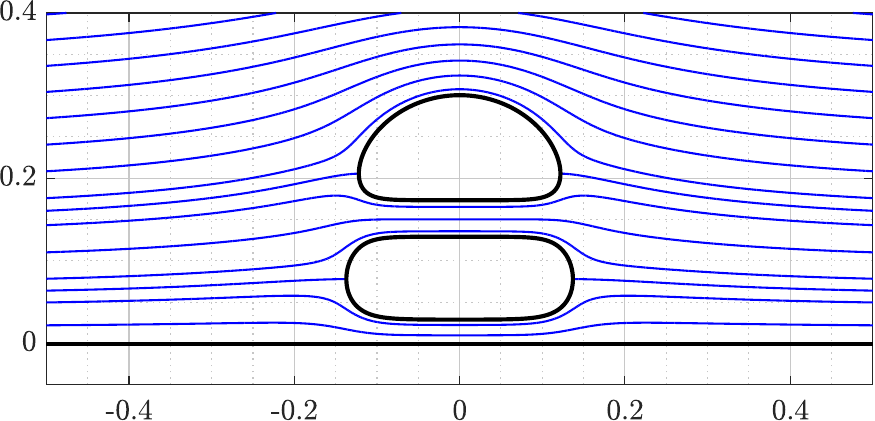}
		\hfill
		\includegraphics[width=0.49\textwidth]{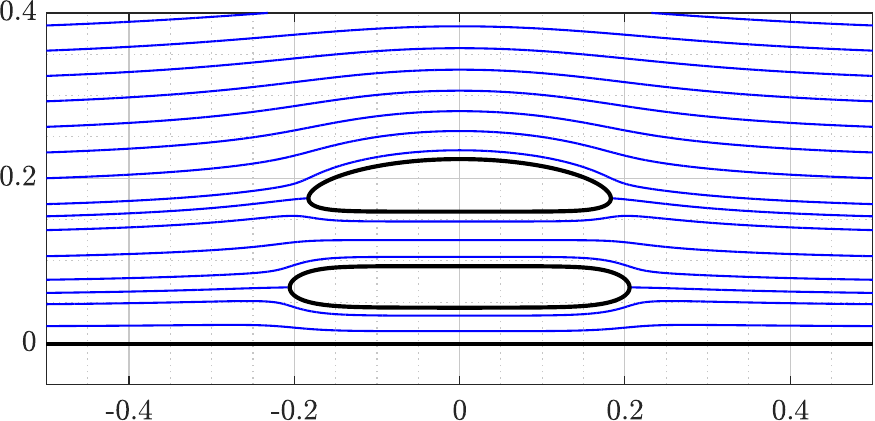}\hfill}
	\caption{The two bubbles in the second case of Example~\ref{ex:2b} with a set of streamlines superposed in the reference frame co-traveling with the bubbles for $U=2$ (left) and $U=4$ (right)}\label{fig:halffig2UW}
\end{figure}

\begin{example}[Two bubbles II]\label{ex:2b}{\rm
Figure~\ref{fig:halffig2U} shows the configurations of two bubbles and several streamlines of the flow field in the co-traveling frame when $U=2$ (left) and $U=4$ (right). 
For this case, the parameters of the circular domain $D_\zeta$ are $r_1=0.0558$, $r_2=0.1003$, $z_1=0.5009\i$, and $z_2=0.3277\i$. These parameters are chosen such that the two bubbles have approximately equal area. Note that in this setup, the bubbles are far from the wall, and their shapes are largely unaffected by the boundary. This leads to a similarity with the configurations observed in the free space case, where no significant external influence (such as proximity to a wall) distorts the boundaries of the bubbles (see Figure~\ref{fig:spac2-lc} above).
		
In contrast, the second case, shown in Figure~\ref{fig:halffig2UW}, illustrates the behavior of the bubbles when they are close to the wall. Here, we observe the interaction between the bubbles and the wall, which causes noticeable distortions in the bubble shapes compared to when they are far from the wall. The parameters of the circular domain $D_\zeta$ for this configuration are $r_1=0.1012$, $r_2=0.1359$, $z_1=-0.5831\i$, and $z_2=-0.8421\i$, chosen so that the two bubbles still have approximately equal area. These configurations reveal the influence of proximity to the boundary on bubble shape, which is more pronounced as the bubbles approach the wall.
}\end{example}

\begin{example}[Two bubbles III]\label{ex:2bb}{\rm
		We again consider two bubbles which are obtained by choosing the parameters of the circular domain $D_\zeta$  such that the two bubbles have approximately equal area. The bubbles are scaled such that their areas are approximately $\pi$.
		Figure~\ref{fig:halffig2Uv} shows the configurations of the two bubbles for several values of translation speed $U$.  The parameters of $D_\zeta$ are $r_1=0.1188288$, $r_2=0.090724$, $z_1=-0.86\i$, and $z_2=-0.62\i$ for the bubbles on the left and $r_1=r_2=0.15$, $z_1=-0.16-0.81\i$, and $z_2=0.16-0.81\i$ for the bubbles on the right.

Comparing with the bubbles in Figure~\ref{fig:spac2}, where the two bubbles always have the same shape, the two bubbles in Figure~\ref{fig:halffig2Uv} (left) no longer have this property owing to the presence of the wall.
The closer bubble to the wall exhibits stronger horizontal elongation illustrating the effect of boundary proximity.		
However, in Figure~\ref{fig:halffig2Uv} (right), we show configurations of two bubbles that have the same shapes. The interactions between the bubbles themselves and with the wall are again evident.  

\begin{figure}[h!]
	\centering
	\centering
	\scalebox{0.43}{\includegraphics[trim=0cm 0.0cm 0cm 0.0cm,clip]{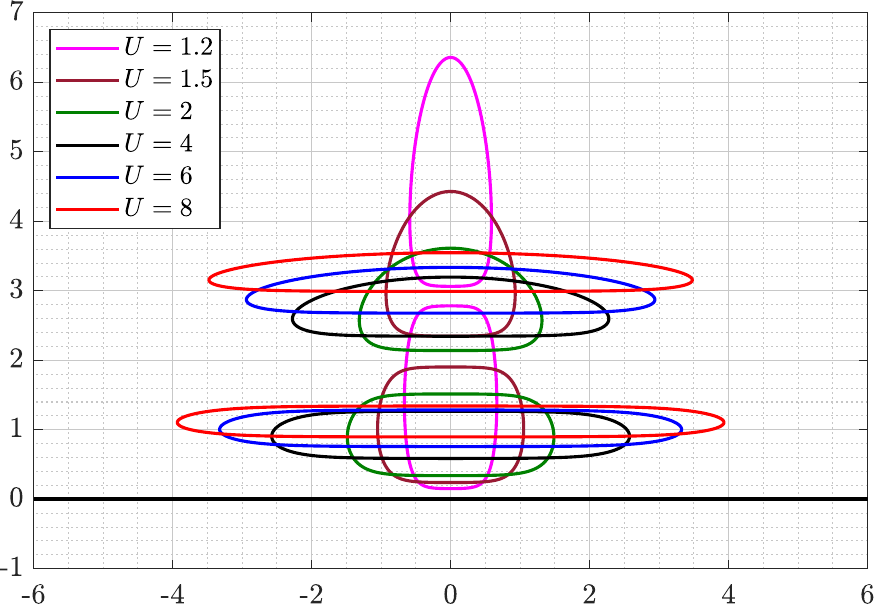}}\hfil
	\scalebox{0.43}{\includegraphics[trim=0cm 0.0cm 0cm 0.0cm,clip]{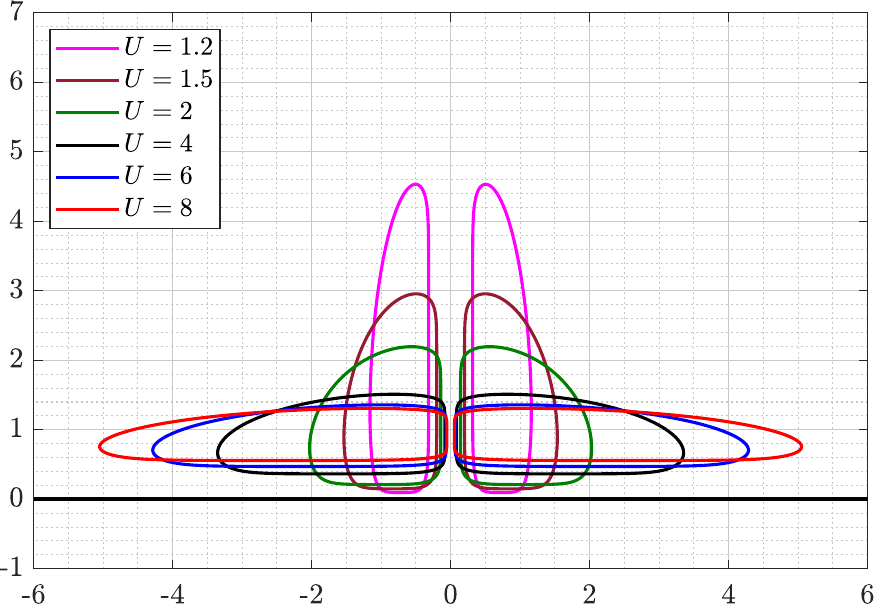}}
	\caption{The configurations of the two bubbles in Example~\ref{ex:2bb} for various values of $U$.}\label{fig:halffig2Uv}
\end{figure}
}\end{example}

\begin{example}[Twenty-five bubbles]\label{ex:25b}{\rm
We present in Figure~\ref{fig:halffigLc25} an asymmetric configuration of $25$ bubbles in the upper half-plane. The figure shows two different configurations with translation speeds $U=2$ (left) and $U=6$ (right). In both cases, a set of streamlines is superimposed in the frame co-traveling with the bubbles, which provides insight into the flow dynamics around the bubble array.		
The configurations highlight how the bubbles interact with the flow field and each other in the upper half-plane. At a higher speed, the interactions between the bubbles themselves become less pronounced. Since the domain is unbounded in the vertical direction, the bubbles furthest from the wall exhibit shapes similar to those seen in free space.

\begin{figure}[h!]
	\centering
	\scalebox{0.43}{\includegraphics[trim=0.0cm 1cm 0.0cm 1.4cm,clip]{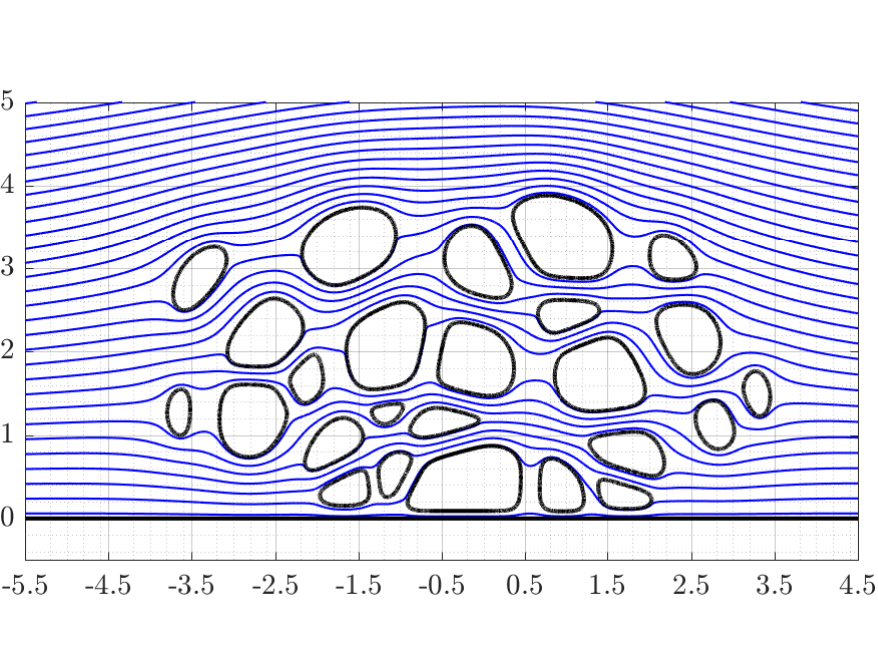}}
	\hfill
	\scalebox{0.43}{\includegraphics[trim=0.0cm 1cm 0.0cm 1.4cm,clip]{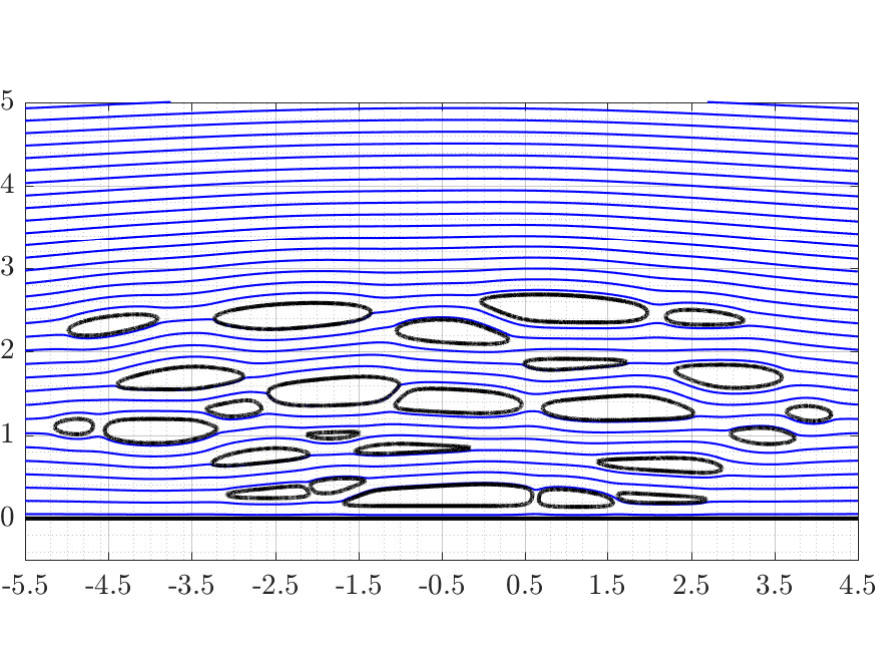}}
	\caption{An example of 25 bubbles in a general asymmetric configuration with a set of streamlines superposed in the reference frame co-traveling with the bubbles when $U=2$ (left) and $U=6$ (right).}\label{fig:halffigLc25}
\end{figure}
}\end{example}

\subsection{Infinite horizontal channel}

In this section, we explore the behavior and configurations of bubbles within an infinite horizontal Hele-Shaw channel. Some of these configurations have been investigated in earlier works, such as Vasconcelos~\cite{vas15}, Green \& Vasconcelos~\cite{GV}, and Tanveer~\cite{tan}.

\begin{figure}[h!]
	\centering
	{\hfill
		\includegraphics[width=0.32\textwidth]{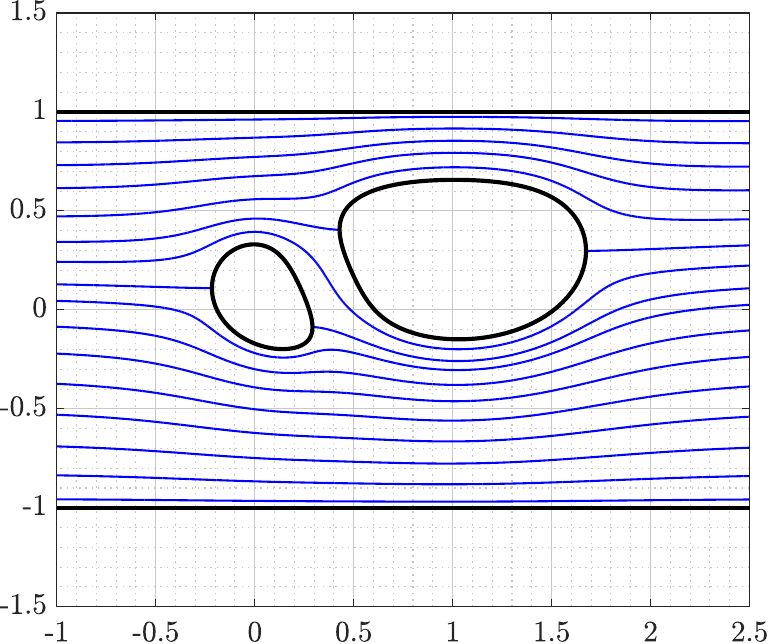}
		\hfill
		\includegraphics[width=0.32\textwidth]{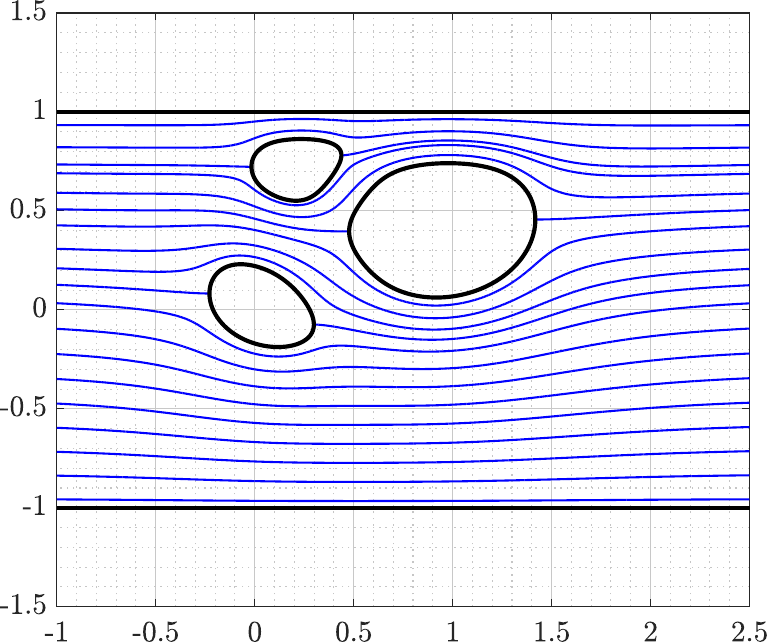}
		\hfill
		\includegraphics[width=0.32\textwidth]{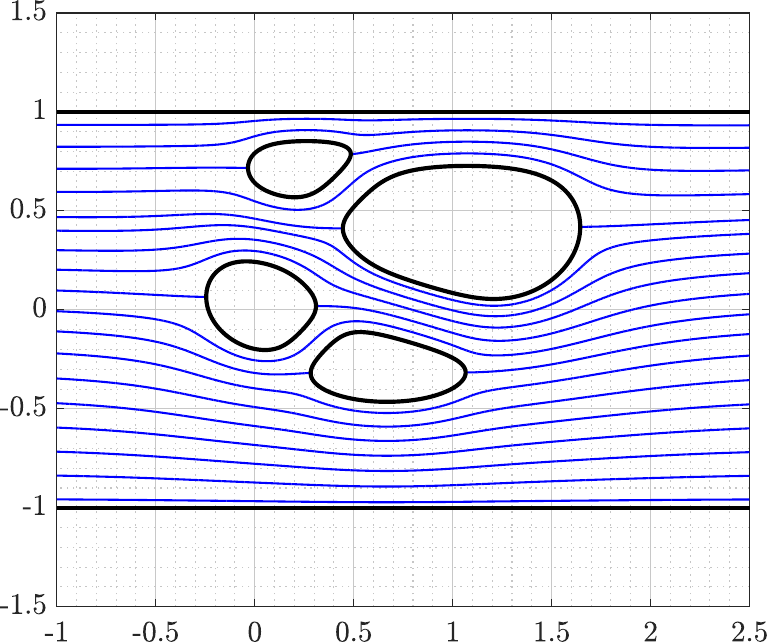}
		\hfill}
	\caption{The configurations of two bubbles (left), three bubbles (middle), and four bubbles (right) in Example~\ref{ex:chan2} with $U=2$.}\label{fig:channel-gv}
\end{figure}

\begin{example}\label{ex:chan2}{\rm
In this example, we recover the bubbles in a Hele-Shaw channel considered in~\cite[Fig. 5--7]{GV}. 
The configurations of two, three, and four bubbles for $U=2$ are shown in Figure~\ref{fig:channel-gv}. We plot streamlines superposed in the frame co-traveling with the bubbles.
The parameters of the circular domains $D_\zeta$ used to generate the bubbles in Figure~\ref{fig:channel-gv} are the same as in~\cite[Fig.~5--7]{GV}. 
For two bubbles, the parameters of the circular domain $D_\zeta$ are: $r_1=0.2$, $r_2=0.25$, $z_1=0.03\i$, and $z_2=0.6+0.15\i$. 
For three bubbles, the parameters of the circular domain $D_\zeta$ are: $r_1=0.18$, $r_2=0.22$, $r_3=0.195$, $z_1=0$, $z_2=0.6+0.23\i$, and $z_3=0.2+0.63\i$. 
For four bubbles, the parameters of the circular domain $D_\zeta$ are: $r_1=0.19$, $r_2=0.235$, $r_3=0.2$, $r_4=0.175$, $z_1=0$, $z_2=0.6+0.23\i$, $z_3=0.2+0.63\i$, and $z_4=0.4-0.25\i$.

	}
\end{example}

\begin{figure}[h!]
\centering
\scalebox{0.43}{\includegraphics[trim=0cm 1.5cm 0cm 1.5cm,clip]{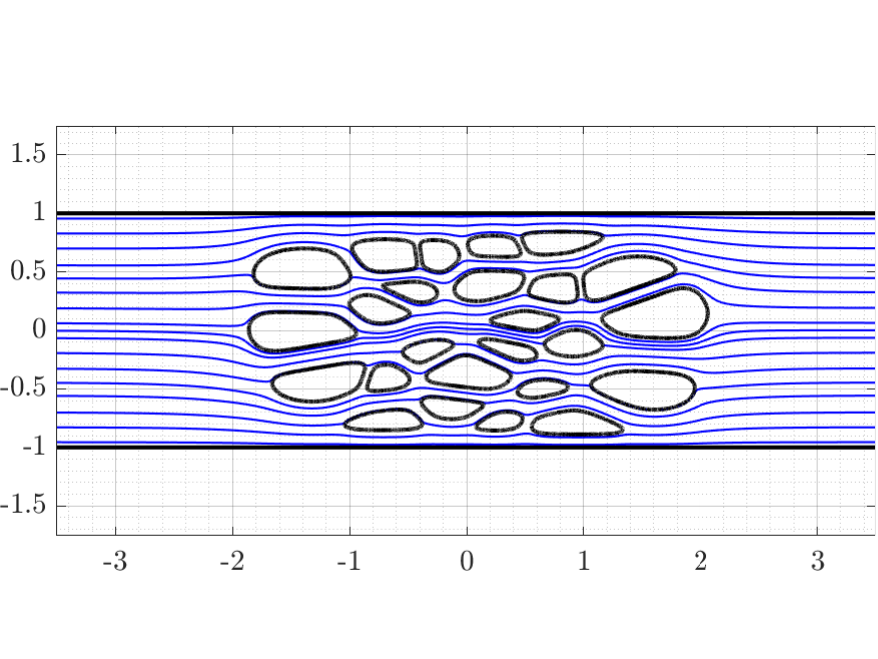}}\hfil
\scalebox{0.43}{\includegraphics[trim=0cm 1.5cm 0cm 1.5cm,clip]{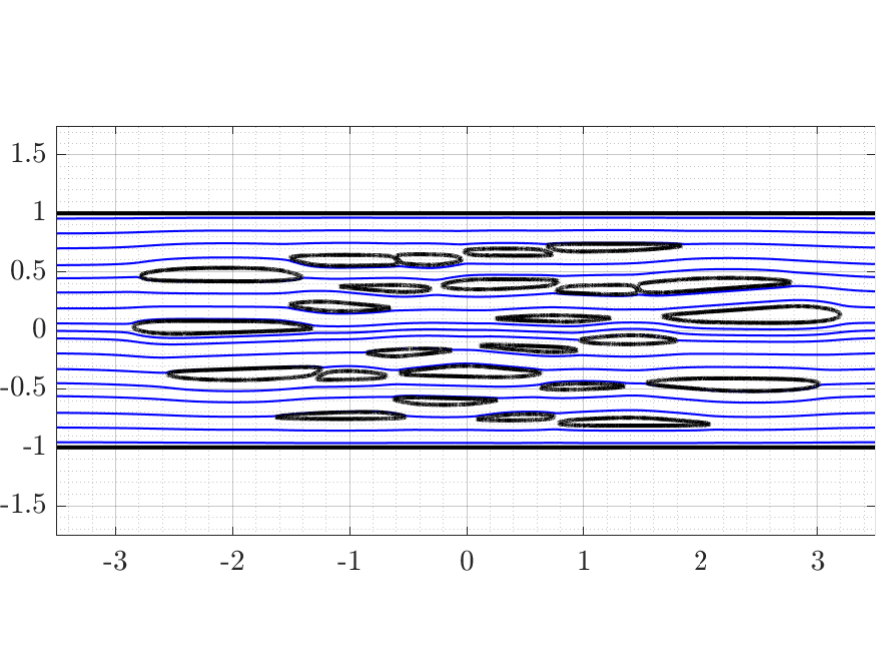}}
\caption{An example of 25 bubbles in a general up-down symmetric configuration in a channel with a set of streamlines superposed in the reference frame co-traveling with the bubbles when $U=2$ (left) and $U=6$ (right).}\label{fig:chanfigLc25}
\end{figure}

\begin{example}[Twenty-five bubbles]\label{ex:chan25}{\rm
Figure~\ref{fig:chanfigLc25} presents an asymmetric configuration of 25 bubbles in a horizontal channel of width 2. This example demonstrates our method's ability to handle complex, multi-bubble scenarios with irregular arrangements. We compare the configurations at two velocities: $U=2$ and $U=6$. We can again observe an increased bubble elongation along the flow direction at higher speed and simplification of streamline patterns as the speed increases. 
These configurations demonstrate how increasing the number of bubbles affects both the flow field and the bubble shapes. The asymmetric positioning of bubbles leads to interesting interactions, both between the bubbles and with the channel walls.}
\end{example}

\begin{example}[Two bubbles I]\label{ex:2bb-ch}{\rm
Figure~\ref{fig:ch-2DU} illustrates the effect of translation speed $U$ on two bubbles of approximately equal area in three different channel configurations: up-down symmetric (top left), left-right symmetric  (top right), and a third different arrangement (bottom). These cases demonstrate how symmetry (or lack thereof) of bubble positioning influences the shapes of the bubbles as $U$ increases.

\begin{figure}[h!]
	\centering
	\scalebox{0.43}{\includegraphics[trim=0cm 1.0cm 0cm 1.0cm,clip]{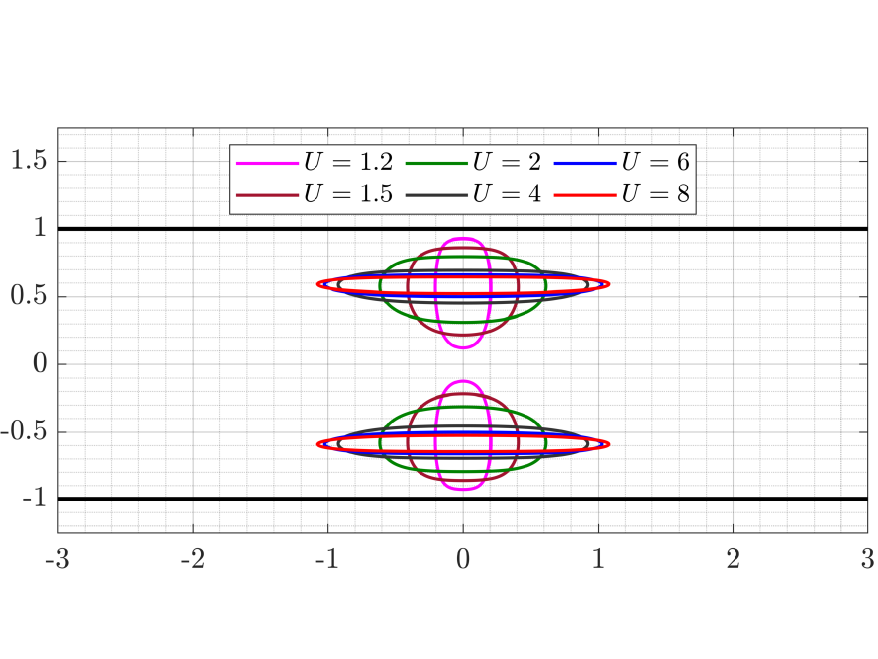}}
	\scalebox{0.43}{\includegraphics[trim=0cm 1.0cm 0cm 1.0cm,clip]{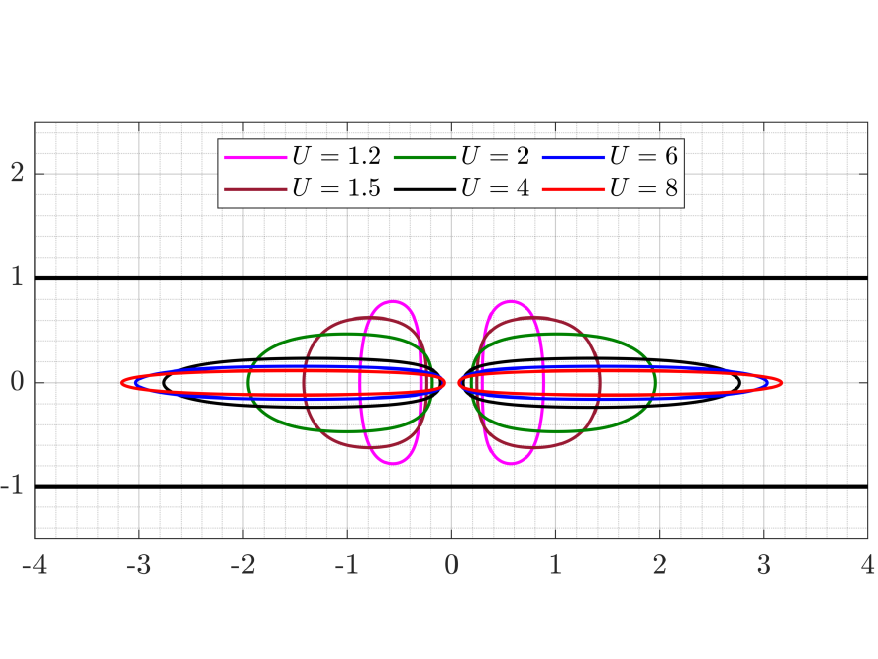}}
	\scalebox{0.43}{\includegraphics[trim=0cm 1.0cm 0cm 1.0cm,clip]{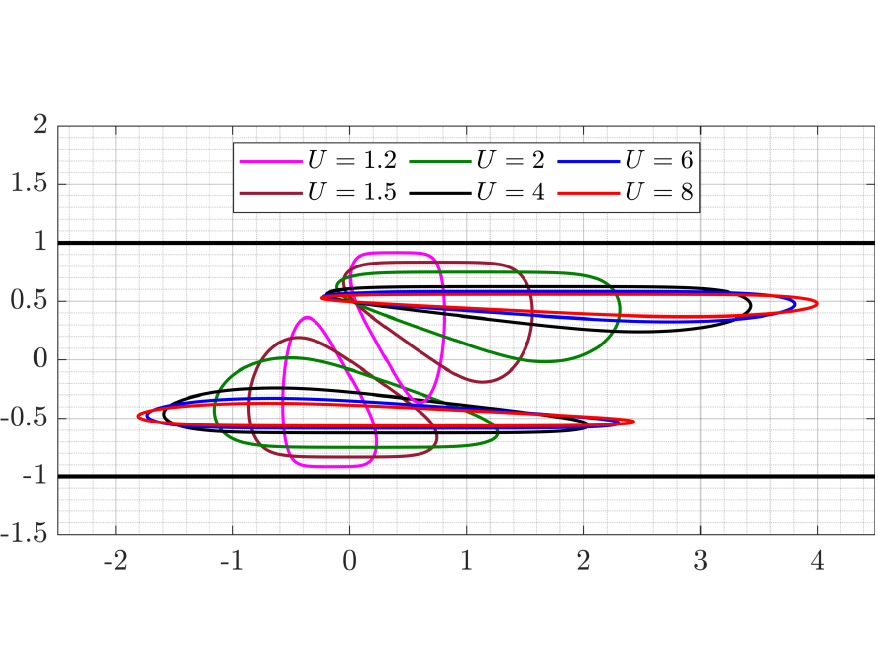}}
	\caption{The configurations of the two bubbles in Example~\ref{ex:2bb-ch} for several values of $U$. The parameters of $D_\zeta$ are:  $r_1=r_2=0.35$ and $z_1=-z_2=0.5\i$ (top left); $r_1=r_2=0.35$ and $z_1=-z_2=0.5$ (top right);  $r_1=r_2=0.42$ and $z_1=-z_2=0.35+0.35\i$  (bottom). }\label{fig:ch-2DU}
\end{figure}

}\end{example}

\begin{example}[Two bubbles II]\label{ex:2bb2-ch}{\rm
In Example~\ref{ex:2bb-ch}, we considered two bubbles which were obtained by fixing the parameters of the circular domain $D_\zeta$ and changing the values of $U$. 
Building on this, we now fix $U=2$ and vary the parameters of the circular domain $D_\zeta$. 
The obtained results are shown in Figure~\ref{fig:ch-2DU2} and reveal how the positions and sizes of the bubbles influence their shapes and interactions at a constant speed.

\begin{figure}[h!]
	\centering
	\scalebox{0.43}{\includegraphics[trim=0cm 1.0cm 0cm 1.5cm,clip]{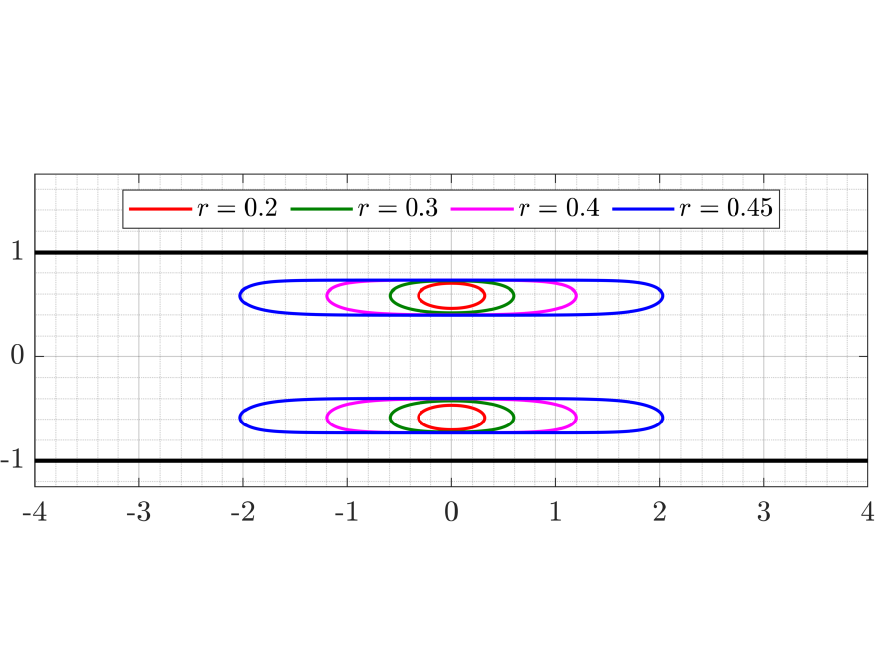}}
	\scalebox{0.43}{\includegraphics[trim=0cm 1.0cm 0cm 1.5cm,clip]{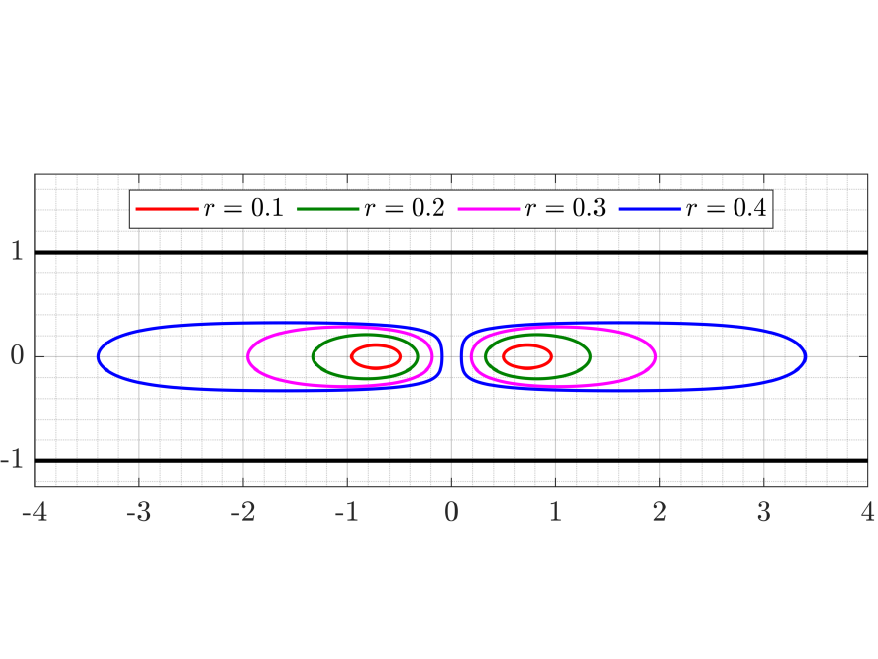}}
	\scalebox{0.43}{\includegraphics[trim=0cm 0.0cm 0cm 0.0cm,clip]{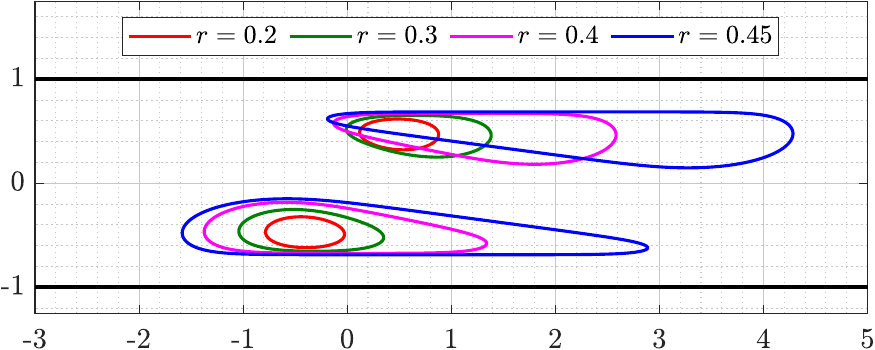}}
	\caption{The configurations of the two bubbles in Example~\ref{ex:2bb2-ch} with $U=2$. The centers of the two inner circles of $D_\zeta$ are:  $z_1=-z_2=0.5\i$ (top left); $z_1=-z_2=0.5$ (top right);  $z_1=-z_2=(1+\i)/2\sqrt{2}$  (bottom). For the three cases, the radii of the two circles are $r_1=r_2=r$ for several values of $r$.}\label{fig:ch-2DU2}
\end{figure}

}\end{example}

\begin{figure}[h!]
	\centering
	\scalebox{0.43}{\includegraphics[trim=0cm 1.0cm 0cm 1.0cm,clip]{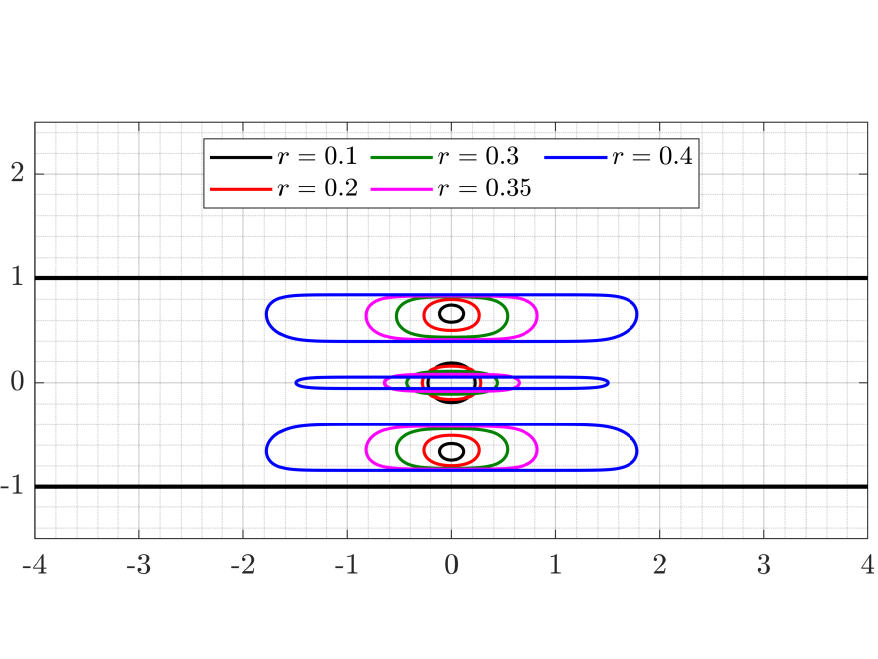}}
	\scalebox{0.43}{\includegraphics[trim=0cm 1.0cm 0cm 1.0cm,clip]{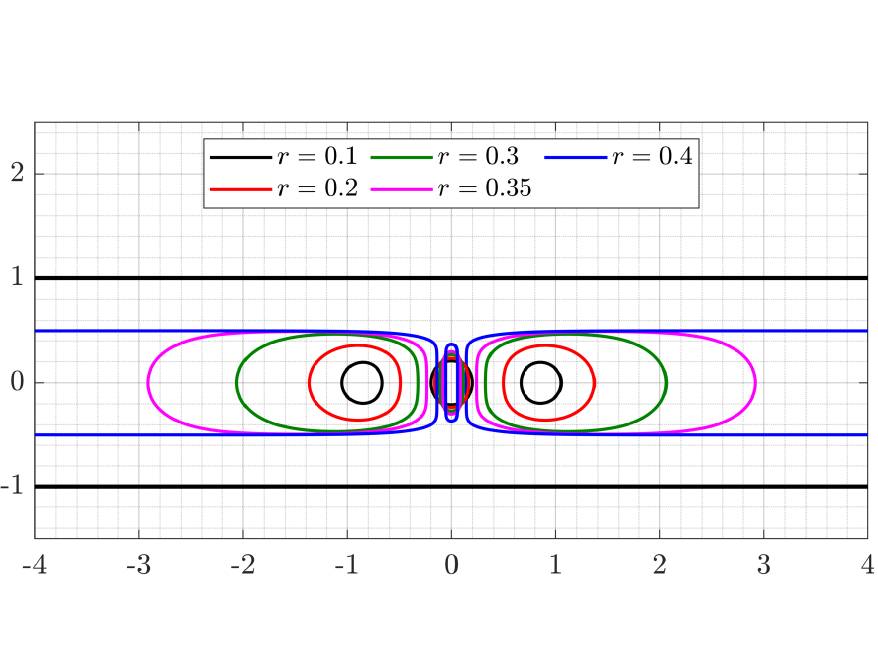}}
	\scalebox{0.43}{\includegraphics[trim=0cm 1.0cm 0cm 1.0cm,clip]{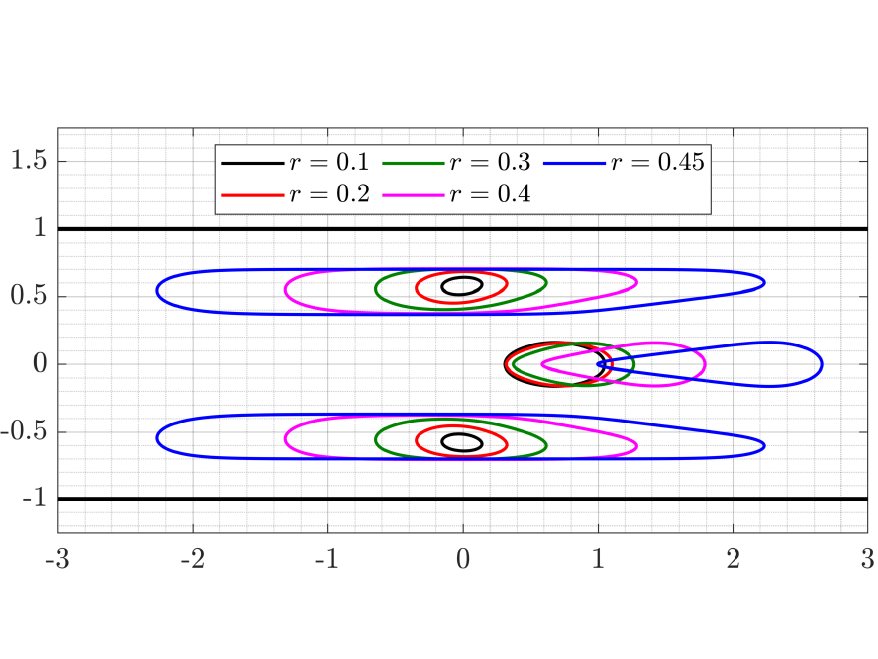}}
	\scalebox{0.43}{\includegraphics[trim=0cm 1.0cm 0cm 1.0cm,clip]{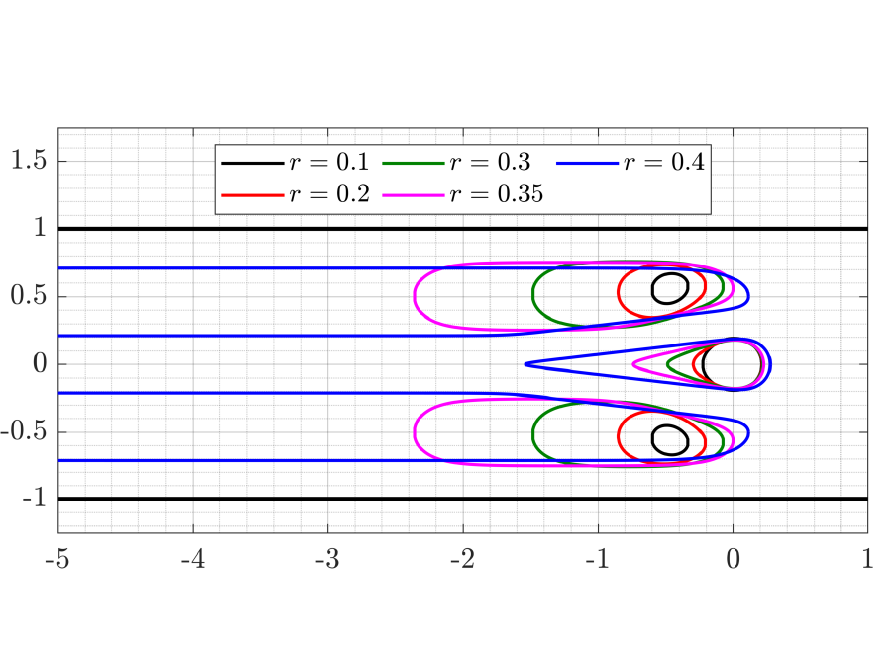}}
	\scalebox{0.43}{\includegraphics[trim=0cm 1.0cm 0cm 1.0cm,clip]{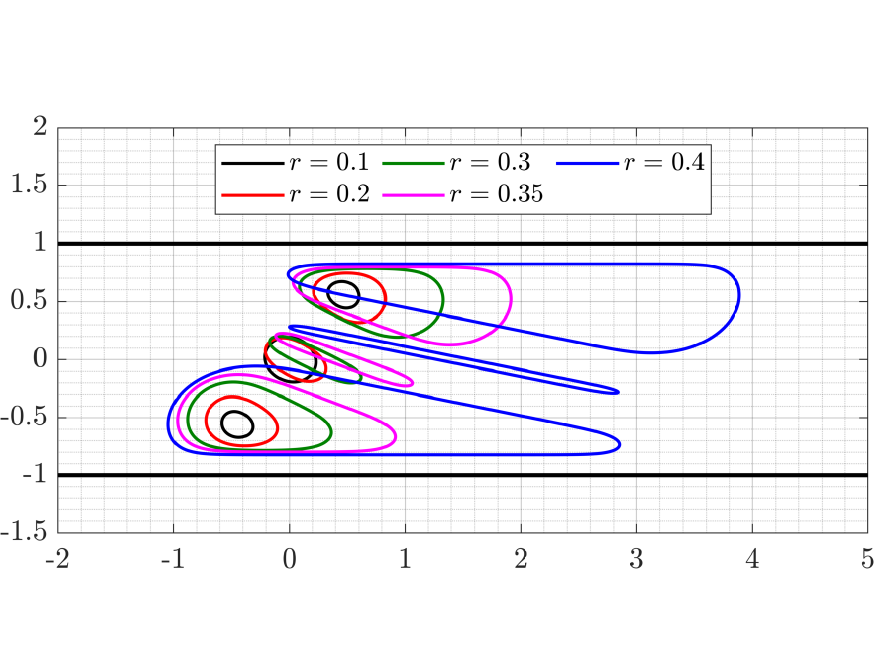}}
	\caption{The configurations of the three bubbles in Example~\ref{ex:3bb-ch} with $U=2$. The centers of the three inner circles of $D_\zeta$ are: 
		$z_1=-z_2=0.58\i$, $z_3=0$ (top left); 
		(b) $z_1=-z_2=0.58$, $z_3=0$ (top right); 
		(c) $z_1=-z_2=0.5\i$, $z_3=0.45$ (middle left); 
		(d) $z_1=\overline{z_2}=-0.41+0.41\i$, $z_3=0$ (middle right); 
		(e) $z_1=-z_2=0.41+0.41\i$, $z_3=0$ (bottom). 
		For all cases, the radii of the three circles are $r_1=r_2=r$ and $r_3=0.16$ for several values of $r$.}\label{fig:channel-3Dr}
\end{figure}

\begin{example}[Three bubbles]\label{ex:3bb-ch}{\rm
We now explore various three-bubble configurations at $U=2$ as shown in Figure~\ref{fig:channel-3Dr} for several values of the parameters of $D_\zeta$.
The configurations of the bubbles are: up-down and left-right symmetric about the coordinate axes (top), up-down symmetric about the horizontal channel center line (middle), and a third different configuration (bottom).

}\end{example}

\begin{figure}[h!]
\centering
\scalebox{0.43}{\includegraphics[trim=0cm 0.0cm 0cm 0.0cm,clip]{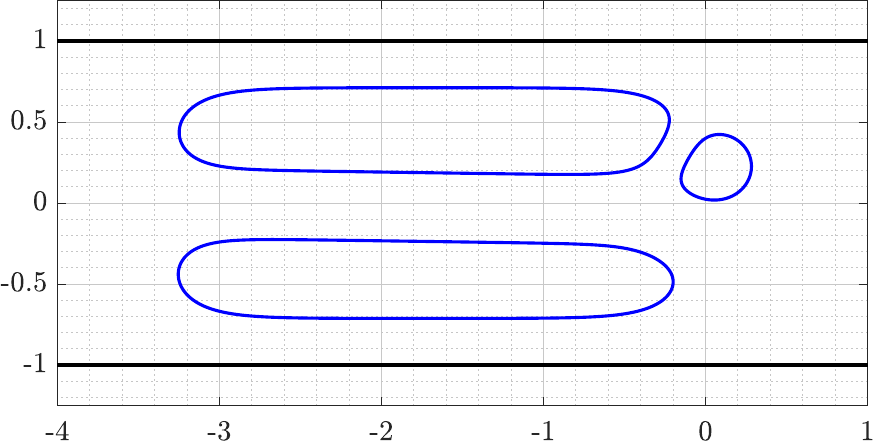}}
\scalebox{0.43}{\includegraphics[trim=0cm 0.0cm 0cm 0.0cm,clip]{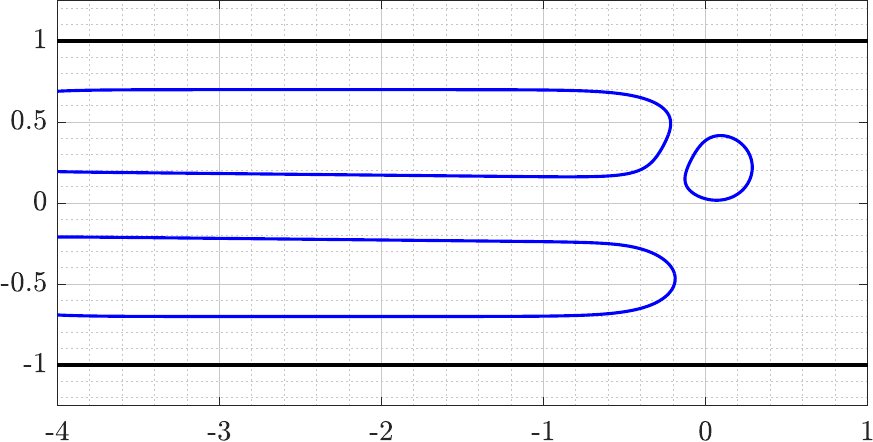}}
\caption{Examples of three bubbles in a horizontal channel whose configurations have been chosen to be qualitatively similar to the configuration of two fingers and one bubble in Figure~9 of Vasconcelos~\cite{vas15} for $U=2$. In both configurations, we have three bubbles, two of which appear to be finger-like.
The parameters of the three inner circles of $D_\zeta$ are: 
$z_1=\overline{z_2}=-0.55+0.33\i$, $z_3=0.05+0.18\i$, $r_1=r_2=0.31$, $r_3=0.33$ (left); 
$z_1=\overline{z_2}=-0.56+0.33\i$, $z_3=0.05+0.18\i$, $r_1=r_2=0.32$, $r_3=0.33$ (right).	
}\label{fig:channel-2F1B}
\end{figure}

\begin{example}{\rm
Our final example presents two special cases of three-bubble configurations designed to qualitatively resemble the two-finger, one-bubble configuration studied by Vasconcelos \cite{vas15}. These configurations, shown in Figure~\ref{fig:channel-2F1B}, illustrate how our method can be applied to study configurations that mimic more complex scenarios, such as the interplay between bubble and finger-like phenomena in Hele-Shaw flows.
}\end{example}
	
\section{Concluding remarks}

\begin{figure}[h!]
	\centering
	\scalebox{0.5}{\includegraphics[trim=0cm 0.0cm 0cm 0.0cm,clip]{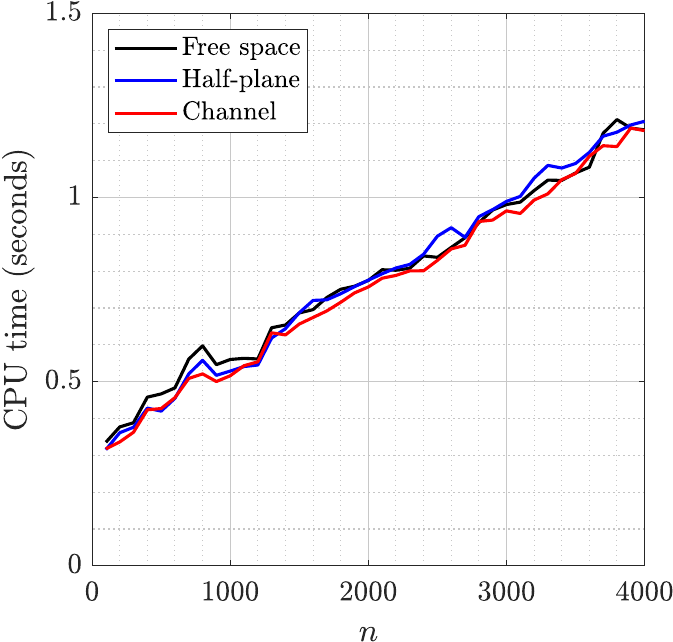}}
	\scalebox{0.5}{\includegraphics[trim=0cm 0.0cm 0cm 0.0cm,clip]{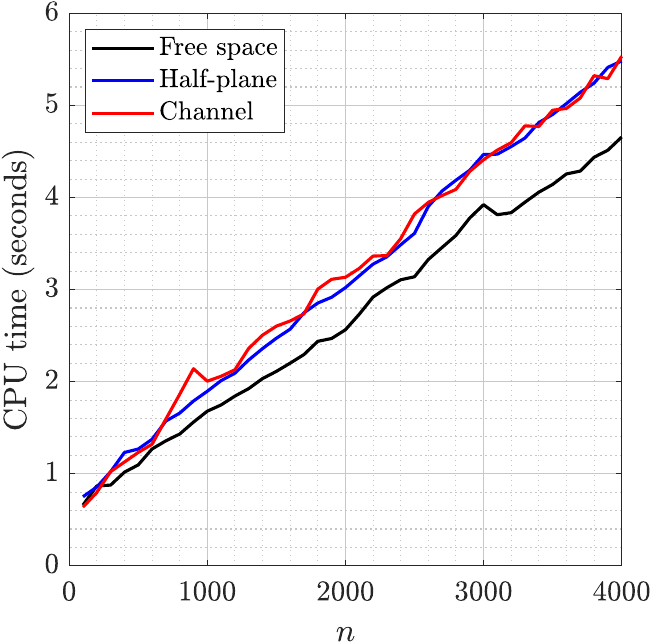}}
	\caption{Run time in seconds as a function of $n$ required to produce: the bubble shapes shown in Figures~\ref{fig:spacefigs}, \ref{fig:halffigs} and~\ref{fig:chanfigs} (left), and the bubble shapes shown in Figures~\ref{fig:spacefigLc25}, \ref{fig:halffigLc25} and~\ref{fig:chanfigLc25} (right).
	}\label{fig:time}
\end{figure}

\begin{figure}[h!]
	\centering
	\scalebox{0.5}{\includegraphics[trim=0cm 0.0cm 0cm 0.0cm,clip]{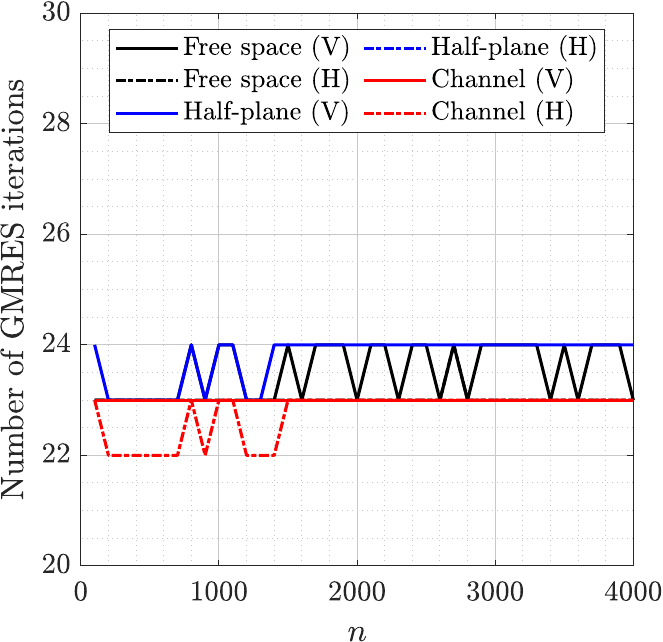}}
	\scalebox{0.5}{\includegraphics[trim=0cm 0.0cm 0cm 0.0cm,clip]{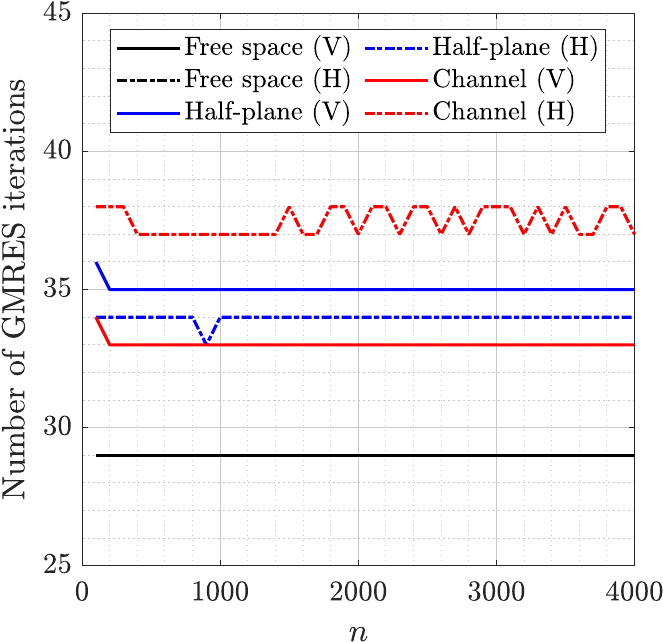}}
	\caption{The number of GMRES iterations as a function of $n$ required to compute the two conformal mappings (one to a canonical domain with horizontal slits (H) and the other to a canonical domain with vertical slits (V)) for the three configurations of interest.}
	\label{fig:gmres}
\end{figure}

Using a fast and accurate boundary integral equation method, combined with conformal mappings, we have presented a unified method for solving the free boundary problems governing multiple bubble boundary shapes in steady Hele-Shaw flow in free space, the upper-half plane, and an infinite horizontal channel. 
Solving each free boundary problem requires the computation of two conformal mappings (one to a canonical domain with horizontal slits and the other to a canonical domain with vertical slits). Computing each conformal mapping requires solving the boundary integral equation~\eqref{eq:ie} as explained briefly in subsection~\ref{sec:numerical}. The computation of the right-hand side of the integral equation requires $\mathcal{O}((m +1)n \ln n)$ operations and each iteration of the GMRES method requires $\mathcal{O}((m +1)n)$ operations, where $m+1$ is the connectivity of the multiply connected domain and $n$ is the number of nodes used in the discretization of each boundary component~\cite{Nas-ETNA}.
To corroborate this, a plot of the CPU time (in seconds) as a function of $n$ required to solve the governing free boundary problem for the bubble shapes is shown in Figure~\ref{fig:time} (left) for the bubbles shown in Figures~\ref{fig:spacefigs}, \ref{fig:halffigs} and~\ref{fig:chanfigs}, and in Figure~\ref{fig:time} (right) for the bubbles shown in Figures~\ref{fig:spacefigLc25}, \ref{fig:halffigLc25} and~\ref{fig:chanfigLc25}. It is clear from these graphs that the CPU time depends almost linearly on $n$. For the $25$ bubble case in Figure~\ref{fig:time} (right), note that $m=24$ for the free space case and $m=25$ for both the upper half-plane and the channel cases which explains why the times for the free space case are lower than the others.
We show also in Figure~\ref{fig:gmres} the number of GMRES iterations required to compute the two conformal mappings used to produce the aforementioned bubble figures. These graphs show that the number of GMRES iterations is almost independent of $n$ (see also~\cite{Nas-ETNA}).

The examples, presented in this paper, illustrate the versatility of our numerical method, which effectively handles a wide range of bubble configurations in a Hele-Shaw channel, from simple two-bubble systems to more complex multi-bubble arrangements. 
As seen in Figures~\ref{fig:spacefigLc25},~\ref{fig:halffigLc25}, \&~\ref{fig:chanfigLc25}, where configurations involving 25 bubbles are depicted, the method remains effective even for highly connected domains. 

In addition to solving problems with complex geometrical configurations, this work provides valuable insights into the dynamics of bubbles in confined geometries, which offers a foundation for future studies on more intricate scenarios in Hele-Shaw flows. However, the free boundary problems considered here have had the major underlying assumption that the surface tension on the bubble boundaries is zero. Recent investigations into selection problems associated with having non-zero surface tension, from two different perspectives~\cite{GLM,MV}, suggest that multi-bubble configurations will introduce further selection mechanisms when surface tension effects are included. This opens additional avenues for exploration in future work.

The MATLAB codes as well as the parameters of the circular domain $D_\zeta$ used to produce the configurations presented in this paper can be found at the link \url{https://github.com/mmsnasser/bubbles}.	

\section*{Acknowledgements} The authors would like to thank an anonymous reviewer for their valuable corrections, comments and suggestions, which greatly improved the presentation of this paper.

\end{document}